\documentstyle[12pt,leqno]{article}
\textheight9in \def\DATE{July 22, 1998}
\newtheorem{theorem}{Theorem}[section]
\newtheorem{definition}[theorem]{Definition}
\newtheorem{corollary}[theorem]{Corollary}

\newtheorem{lemma}[theorem]{Lemma}

\newtheorem{proposition}[theorem]{Proposition}

\catcode`\@=11

\def\ps@myheadings{\let\@mkboth\@gobbletwo
\def\@oddhead{\ifnum\count0=1 \hfill\else
\rightmark \hfil \rm\thepage\fi}%
\def\@oddfoot{\ifnum\count0=1 \hfill \rm 1 \hfill \else
\hfill\fi}
\def\@evenhead%
{\rm\leftmark\hfil\rm\thepage}%
\def\@evenfoot{}\def\sectionmark##1{}
\def\subsectionmark##1{}}

\def\@begintheorem#1#2{\it \trivlist \item[\hskip
 \labelsep{\bf #1\ #2.}]}
\def\@opargbegintheorem#1#2#3{\it \trivlist\item[\hskip%
 \labelsep{\bf #1\ #2.\ (#3)}]}
\def\@endtheorem{\endtrivlist}

\def\@listI{\leftmargin\leftmargini \parsep 1pt plus 2.5pt
 minus 1pt\topsep 5pt plus 2pt minus 3pt%
 \itemsep 0pt plus 2.5pt minus 1pt}
\let\@listi\@listI
\@listi

\def\@sect#1#2#3#4#5#6[#7]#8{\ifnum #2>\c@secnumdepth%
 \def \@svsec {}\else \refstepcounter {#1}\edef \@svsec%
 {\csname the#1\endcsname. \hskip .1em }\fi \@tempskipa%
 #5\relax \ifdim \@tempskipa >\z@ \begingroup #6\relax%
 \@hangfrom {\hskip #3\relax \@svsec }{\interlinepenalty%
 \@M #8.\par }\endgroup \csname #1mark\endcsname {#7}%
 \addcontentsline {toc}{#1}{\ifnum #2>\c@secnumdepth%
 \else \protect \numberline {\csname the#1\endcsname. }%
 \fi #7}\else \def \@svsechd {#6\hskip #3\@svsec #8.%
 \csname #1mark\endcsname {#7}\addcontentsline {toc}{#1}%
 {\ifnum #2>\c@secnumdepth \else \protect \numberline%
 {\csname the#1\endcsname. }\fi #7}}\fi \@xsect {#5}}

\def\section{\@startsection {section}{1}{\z@ }%
 {-3.5ex plus -1ex minus -.2ex}{2.3ex plus .2ex}{\bf }}

\def\thebibliography#1{%
 \section *{References.\@mkboth {REFERENCES}{REFERENCES}}%
 \list {[\arabic {enumi}]}{\settowidth \labelwidth {[#1]}%
 \leftmargin \labelwidth \advance \leftmargin \labelsep %
 \usecounter {enumi}} \def \newblock %
 {\hskip .11em plus .33em minus -.07em} \sloppy \clubpenalty 4000%
 \widowpenalty 4000 \sfcode`\.=1000\relax}

\def\@maketitle{%
 \newpage \null \vskip 2em
 \begin{center}
{\Large\bf \@title \par }
 \vskip 1.5em
 {\large \lineskip .5em
 \begin {tabular}[t]{c}\@author
 \end{tabular}\par}
 \end{center}
  \vskip .8em}

\def\abstract{%
\if@twocolumn \section *{Abstract}
 \else \small\quotation\noindent{\bf Abstract.}\fi}

\catcode`\@=13

\newcommand{\square}[8]{
\setlength{\unitlength}{.7cm}
\begin{picture}(5,3.1)
\thicklines

\put(0,3){\makebox(0,0){$#1$}}
\put(5,3){\makebox(0,0){$#2$}}
\put(0,0){\makebox(0,0){$#3$}}
\put(5,0){\makebox(0,0){$#4$}}

\put(-.5,1.5){\makebox(0,0)[r]{$#6$}}
\put(5.5,1.5){\makebox(0,0)[l]{$#7$}}
\put(2.5,0.5){\makebox(0,0)[b]{$#8$}}
\put(2.5,3.5){\makebox(0,0)[b]{$#5$}}

\put(1,0){\vector(1,0){3}}
\put(1,3){\vector(1,0){3}}
\put(0,2.5){\vector(0,-1){2}}
\put(5,2.5){\vector(0,-1){2}}
\end{picture}
}

\def\znamenko#1{{(-1)^{(#1)}\cdot}}
\def\EE#1#2{{E^{#1,#2}}}
\def\pa{\partial}
\def\ccdots{{\hskip-1mm \cdot\!\cdot\!\cdot \hskip-1mm}}
\def\x{{\bf x}}
\def\bfr{{\bf R}}
\def\unsh{{\it unsh}}
\def\DO{{\rm DO}}
\def\ldo{{\rm LDO}}
\def\dend{{\rm DEnd}}
\def\ot{{\otimes}}
\def\B{{\cal B}}

\def\boxed#1{\fbox{\rm #1}\hskip1mm}

\def\u{{\bf u}}
\def\y{{\bf y}}

\def\doubless#1#2{{
\def\arraystretch{.5}
\begin{array}{c}
\mbox{\scriptsize $\scriptstyle #1$}
\\
\mbox{\scriptsize $\scriptstyle #2$}
\end{array}\def\arraystretch{1}
}}

\def\prada#1#2{{#1+\cdots+#2}}
\def\orada#1#2{{#1\otimes\cdots\otimes#2}}

\def\Rada#1#2#3{#1_{#2},\dots,#1_{#3}}

\def\Om#1{{\Omega_{#1}(J^\infty E)}}
\def\eps{{\epsilon}}
\def\sgn#1{{(-1)^{#1}}}
\def\dr#1#2{{\Omega^{#1}({\bf R}^{#2})}}
\def\el#1{{\Omega^{#1,0}(J^\infty E)}}
\def\ci#1{{C^{\infty}({\bf R}^{#1})}}
\def\om#1{{\Omega^{#1}(J^\infty E)}}
\def\rada#1#2{#1,\ldots,#2}
\def\prada#1#2{#1+\cdots+#2}

\def\missing#1{{\widehat{#1}}}
\def\loc{{\rm Loc}}
\def\coll#1{{\{#1(n)\}_{n\geq 1}}}

\def\cases#1#2#3#4{
                  \def\arraystretch{1.4}
                         \left\{
                            \begin{array}{ll}
                             #1,\ &\mbox{#2}
                             \\
                             #3,\ &\mbox{#4}
                          \end{array}
                         \right.
                   \def\arraystretch{1}
}
\def\tricases#1#2#3#4#5#6{
                  \def\arraystretch{1.4}
                         \left\{
                            \begin{array}{ll}
                             #1,\ &\mbox{#2}
                             \\
                             #3,\ &\mbox{#4}
                             \\
                             #5,\ &\mbox{#6}
                          \end{array}
                         \right.
                   \def\arraystretch{1}
}

\def\fr#1#2{\ifinner {#1}/{#2} \else \frac{#1}{#2}\fi}
\def\bleft{\ifinner ( \else \left(\rule{0pt}{14pt} \right. \fi}
\def\bright{\ifinner ) \else \left.\rule{0pt}{14pt} \right)\fi}
\def\bleftsq{\ifinner [ \else \left[\rule{0pt}{14pt} \right. \fi}
\def\brightsq{\ifinner ] \else \left.\rule{0pt}{14pt} \right]\fi}

\def\qed{\hspace*{\fill}
\mbox{\hphantom{mm}\rule{0.25cm}{0.25cm}}\\}

\begin{document}

\bibliographystyle{plain}
\baselineskip 16pt plus 2pt minus 1 pt

\title{Differential Operator Endomorphisms of an 
Euler-Lagrange Complex}

\author{Martin Markl%
\thanks{This author was supported by the
grant GA~\v CR 201/96/0310}\hphantom{.}
and Steve Shnider}

\maketitle

\begin{abstract}
The  main  results of our paper
deal with the lifting problem for multilinear differential 
operators between complexes of horizontal de~Rham forms on the
infinite jet bundle. 
We answer the question 
when does an $n$-multilinear differential operator from 
the space of $(N,0)$-forms (where $N$ is the dimension of the base) 
to the space of $(N-s,0)$-forms
allow an $n$-multilinear extension of degree $(-s,0)$ defined on the
whole horizontal de~Rham complex. To study this problem we define
a differential graded operad $\dend_*$ of multilinear differential
endomorphisms, which we prove (Theorem~\ref{Mixture}) to be acyclic 
in positive degrees (negative mapping degrees)
and describe the cohomology group in degree zero in terms
of the characteristic (Definition \ref{defchar}).
Corollary~\ref{Husicka} uses this result to solve
the lifting problem. An important application to mathematical physics
is the proof of  existence of a strongly homotopy Lie algebra
structure extending 
a Lie bracket on the space of functionals (Theorem~\ref{BFLS}).

\noindent
The results of the paper were announced in~\cite{markl-shnider:lifting}.
\end{abstract}

\halign{
\vtop{\parindent=0pt\hsize=10em\strut#\strut\hfill}&
\vtop{\parindent=0pt\hsize=5.5in\strut#\strut}\cr
{\bf Plan of the paper:} & \ref{introduction}.~Introduction\cr
& \ref{Mi}.~First toys\cr
& \ref{ra}.~More serious models\cr
& \ref{main_results}.~Main results\cr
& \ref{proofs}.~Proofs\cr
& \ref{applications}.~Applications\cr
}

\noindent 
{\bf Classification:}
53B50, 58A12, 58A20

\noindent
{\bf Keywords:}
partial differential operator, Euler-Lagrange complex, local
functional

\section{Introduction}
\label{introduction}
Our interest in the subject began when we read the paper
\cite{barnich-fulp-lada-stasheff:preprint} in which the authors
construct a strong homotopy Lie algebra extending the Poisson bracket on
local functionals. The horizontal de~Rham complex on the infinite
jet bundle
(one row of the variational bicomplex) can be augmented over the
space of local functionals by a map defined as integration of
the pull-back of an $(N,0)$-form by a section of the jet bundle,
see Definition~\ref{localfunctional}.
This  defines a projective complex which, after passing to a quotient
by the space of constants in the bidegree $(0,0)$ term, gives
 a resolution of the space of local 
functionals. One can define an extension of the Poisson bracket
by applying standard techniques for resolutions (modulo constants).
The problem with this approach is that the extension is done
value by value with no control of the type of
operator (differential, continuous, etc.) being defined;
when the Poisson bracket is given by a differential operator,
the higher brackets may not be,  and, in fact, need
 not have any particular regularity properties.

We use  a different approach for which all extensions will belong
to a natural class of differential operators, ``local differential
operators," which is the one to which the
Poisson brackets are usually assumed to belong and includes the 
Euler-Lagrange operator and the total horizontal derivatives,
see~\ref{sa} for the definition and explanation of the terminology.

The precise statements of
our results are both technically and notationally complicated,
and therefore, we have decided to motivate the reader by some `toy
models' -- the case of linear (i.e.~not multilinear) operators on
the zero-dimensional bundle ${\bf R}^N \to {\bf R}^N$,
where the forms on the jet
bundle are, of course, ordinary de~Rham forms on ${\bf R}^N$.
In Section~\ref{ra} we add vertical variables, but still remain in
the linear case. The multilinear situation is introduced in
Section~\ref{main_results}.

\noindent
{\bf Acknowledgment.}
The authors would like to express their thanks to Jim Stasheff for
reading the manuscript and many useful comments and suggestions.

\section{First toys}
\label{Mi}

Let us consider the de~Rham complex $\Omega^0({\bf R}) 
\stackrel {d}{\longrightarrow} \Omega^1({\bf R})$
on the one-dimensional Euclidean space ${\bf R}$.
The space, $\DO({\bf R})$,  of linear differential 
operators on ${\bf R}$, consists of maps
\begin{eqnarray*}
A&=&\sum_{i\geq 0}a_i\left(\frac{d}{dx}\right)^i:
C^{\infty}({\bf R})\rightarrow C^{\infty}({\bf R}),
\\
A(f) &=& \sum_{i\geq 0}a_i\fr{d^i f}{dx^i}, 
\mbox { for $f \in C^{\infty}({\bf R})$,}
\end{eqnarray*}
where $a_i = a_i(x)$, $i\geq 0$,
is a sequence of smooth functions such that
$a_i=0$ for $i$ sufficiently large.
Any differential operator $A$ has a natural extension to
the space of one forms, 
$A_1:\Omega^1({\bf R})\rightarrow \Omega^1({\bf R})$, given by 
\begin{equation}
\label{Pejsek_s_Kocickou}
A_1(f dx) := A(f)dx \in  \Omega^1({\bf R}).
\end{equation}
We are looking for a  differential operator $A_0 :
\Omega^0({\bf R})
\to  \Omega^0({\bf R})$,
$A_0(g) = \sum_{j\geq 0}b_j \fr{d^j g}{dx^j}$,
which lifts $A_1$ in the sense that $d A_0 =
A_1 d$ or, diagrammatically,
\begin{eqnarray}
\label{enhanced}
\\
\nonumber
\square{\Omega^0({\bf R})}{\Omega^1({\bf R})}{\Omega^0({\bf R})}%
{\Omega^1({\bf R})}d{A_0}{A_1}d
\end{eqnarray}

\noindent
The condition $d A_0 = A_1d$ can be expanded into the system
\begin{eqnarray}
\label{2}
0 &=& \frac{db_0}{dx}
\\
\nonumber
a_0 &=& b_0 + \frac{db_1}{dx}
\\ \nonumber
&\vdots&
\\ \nonumber
a_{k-1} &=& b_{k-1} + \frac{db_k}{dx},\ k\geq 1.
\end{eqnarray}

\noindent
Assuming that $A_1$ has order $n$ and $A_0$ has finite order,
we get the
solution
\begin{eqnarray}
\nonumber b_n &=& a_n
\\ \nonumber
&\vdots&
\\ \nonumber
b_{n-k}&=& a_{n-k} + \sum_{1\leq j\leq k}(-1)^j \frac{d^j
a_{n-k+j}}{dx^j},
\ 1\leq k\leq n.
\\ 
\label{2n}
0 &=& \sum_{1\leq j\leq n+1}(-1)^j \frac{d^j a_{j-1}}{dx^j}.
\end{eqnarray}

Equation (\ref{2n}) imposes the only restriction on
the operators, i.e.~on the coefficients $a_i$ and $b_j$,
for which a lifting exists in the context of our toy
model. Appropriate
generalizations of this condition will reappear in all subsequent
examples.
The following definitions will allow a precise formulation of
necessary
and sufficient conditions for a lifting.

\begin{definition}The {\em formal adjoint} of the differential
operator
\[
A=\sum_{i\geq 0}a_i\left(\frac{d}{dx}\right)^i
\]
is defined as the differential operator
\[
A^+:=\sum_{i\geq 0}(-1)^i\left(\frac{d}{dx}\right)^i\circ a_i,
\]
where $(\fr{d}{dx})^i\circ a_i$ is the  composition of 
the operator of the $i$-th derivative 
with multiplication by $a_i$.
The {\em characteristic\/} of $A$ is the function 
$\chi:\DO({\bf R})\longrightarrow \ci{}$ given by
\[
\chi(A):=A^+(1)=\sum_{i\geq 0}(-1)^i\fr{d^{i} a_i}{dx^i}\in \ci{},
\]
 where $1$ denotes the constant function.
\end{definition}

\noindent
As a consequence of the equation $ (AB)^+=B^+A^+$
for adjoints, the characteristic satisfies
\begin{equation}
\label{zha}
\chi(A\circ B) = B^+(\chi(A)),
\end{equation}
see~\cite{zharinov:92}.
As an immediate consequence of~(\ref{zha}), we obtain
\begin{equation}\label{Oslicek}
\chi(A \circ \frac d{dx}) = -\frac d{dx}\chi(A)\
\mbox{ and }\
\chi(\frac d{dx} \circ A) = 0.
\end{equation}
An important, though obvious,
property is that $\chi$ is a projector, $\chi^2 = \chi$.

\begin{proposition}
\label{sysrq}
For an arbitrary differential operator
$A:\dr 1{} \to \dr 1{}\cong \ci{}$,
there exists a differential operator $\tilde A: \dr 1{}
\to \dr 0{}\cong \ci{}$ such that
\begin{equation}
\label{petr}
A = d\circ \tilde A + \chi(A).
\end{equation}
\end{proposition}

\noindent
{\bf Proof.}
Clearly it is enough to prove the proposition for  $A(f)=
(a_n\fr{d^nf}{dx^n})dx$, for an arbitrary $n\geq 0$. For such $A$,
(\ref{petr}) is
satisfied with
\[
\tilde A(f): = \sum_{0\leq i< n}(-1)^i
\fr{d^ia_n}{d x^i}
\cdot\fr{d^{n-i-1}f}{dx^{n-i-1}}.
\]
\vskip-1cm
\qed

\begin{corollary}
\label{rollex}
In the situation of~(\ref{enhanced}),
the lift $A_0$ of the operator $A_1$ exists if and only if
\[
\chi(A_1) = \mbox{constant}.
\]
\end{corollary}

\noindent
{\bf Proof.}
If $A_1 d=dA_0$ for some
$A_0$, applying  Proposition~\ref{sysrq} to $A=A_1$ and composing
on the right with $d$ gives
\[
d\circ A_0=A_1 \circ d = (d\circ\tilde A + \chi(A_1))\circ
d=d\circ\tilde A
\circ d + \chi(A_1)\circ d.
\]
{}From the last equation, the projector property of $\chi$
and~(\ref{Oslicek}), we have
\[ -\frac{d}{dx}
(\chi(A_1))=\chi(\chi(A_1)\circ d) =\chi(d\circ(A_0-\tilde
A\circ  d))=0,
\]
so $\chi(A_1)=\mbox{constant}$. On the other hand, if the
assumption of
the corollary  is true, we
can put $A_0:= \tilde A\circ d + \chi(A_1)$ and
$d\circ A_0=A_1 \circ d$. \qed

\noindent
Let us move on to a higher-dimensional version of the above
situation. The following notation is standard:
 \begin{equation}
\label{mx}
\left(\frac{\pa}{\pa \x}\right)^I :=
\left(\frac\pa{\pa x^1}\right)^{i_1}\cdots
\left(\frac{\pa}{\pa x^N}\right)^{i_N},
\end{equation}
where $I =(\rada {i_1}{i_N})$,
$\rada{i_1}{i_N} \geq 0$.
\begin{definition}
\label{genius}
Given a linear (partial) differential operator
%, that is, a linear map $A:\ci N \to \ci N$
\begin{equation}
\label{kos}
A = \sum_Ia_I\left(\frac{\pa}{\pa \x}\right)^I 
%f,\f\in \ci N,
\end{equation}
where $a_I = a_I(\x)$ are smooth functions on ${\bf R}^N$
%of\ $\x = (\rada {x^1}{x^N})\in{\bf R}^N$
and
\begin{equation}
\label{cr}
\mbox{$a_I \not= 0$ only for finitely many indices $I$,}
\end{equation}
we define the {\em characteristic\/} to be the function
\begin{equation}
\label{zapisnicek}
\chi(A) := \sum_I (-1)^I \left(\frac{\pa}{\pa \x}\right)^I
a_I \in \ci N,
\end{equation}
where $\sgn I := \sgn{\prada {i_1}{i_N}}$.
\end{definition}

There is a natural definition of a linear partial
differential operator on the space of de~Rham forms. 
The space $\dr kN$ is a free  $\ci N$-module
with basis 
\[
\{(d\x)^\epsilon := (dx^1)^{\epsilon^1}\wedge \cdots
 \wedge(dx^N)^{\epsilon^N}
;\
\epsilon = (\rada {\epsilon^1}{\epsilon^N}),\
\rada{\epsilon^1}{\epsilon^N} \in \{0,1\},\
|\epsilon|:=\sum \epsilon^i = k\}.
\]

\begin{definition}
\label{bednicky}
A linear map $A: \dr kN \to \dr lN$ is a differential operator if
in the expansion
\[
\sum_\epsilon A(f_\epsilon (d\x)^\eps) = \sum_{\epsilon, \delta}
A^\epsilon_\delta (f_\epsilon)(d\x)^\delta,
\]
the `matrix elements' $A^\epsilon_\delta$ are differential
operators in the sense of Definition~\ref{genius}.
\end{definition}

\noindent
Let us consider the complex of de~Rham
forms on ${\bf R}^N$, $N\geq 1$:
\[
0\longrightarrow {\bf R} \longrightarrow
\dr 0N \stackrel {d}{\longrightarrow}
\dr 1N \stackrel {d}{\longrightarrow} \cdots \stackrel
{d}{\longrightarrow} \dr {N-1}N \stackrel {d}{\longrightarrow}
\dr NN
\longrightarrow 0
\]
and a differential operator $A :\dr NN \to \dr NN\cong \ci N $.
The following statement  generalizes
Proposition~\ref{sysrq} and is proved by an easy induction on the number of
variables. 

\begin{proposition}
\label{3}
For any differential operator $A:\dr NN \to \dr NN$ as defined
above, there exists a differential operator $\tilde A: \dr NN
\to \dr
{N-1}N$ such that
\[
A = d \circ \tilde A + \chi(A).
\]
\end{proposition}

\begin{corollary}
A differential operator $A_N : \dr NN \to \dr NN$
can be lifted to a sequence of differential
operators $\{A_s : \dr sN \to \dr sN\}_{0\leq s\leq N}$ such that
$dA_s= A_{s+1}d$ if and only if
\[
\chi(A_N) = \mbox{constant}.
\]
In this case, the lift can be chosen in such a way that
\begin{equation}
\label{favorit}
A_s = \chi(A_N) \mbox{ (multiplication by $\chi(A_N)$)},
\end{equation}
for $0\leq s\leq N-2$.
\end{corollary}

\noindent
{\bf Proof.} Our situation is described by the following diagram:
\begin{center}
{% Picture saved by xtexcad 2.4
\unitlength=1.000000pt
\begin{picture}(200.00,95.00)(0.00,0.00)
\put(194.50,70.50){\vector(0,-1){51.00}}
\put(115.00,50.00){\line(0,-1){10.00}}
\put(115.00,70.00){\line(0,-1){10.00}}
\put(0.00,50.00){\line(0,-1){10.00}}
\put(0.00,60.00){\line(0,1){10.00}}
\put(35.00,10.00){\makebox(0.00,0.00){$d$}}
\put(35.00,90.00){\makebox(0.00,0.00){$d$}}
\put(12.50,40.00){\makebox(0.00,0.00)[l]{$A_0$}}
\put(0.00,30.00){\vector(0,-1){10.00}}
\put(0.00,0.00){\makebox(0.00,0.00){$\dr 0N$}}
\put(0.00,80.00){\makebox(0.00,0.00){$\dr 0N$}}
\put(25.00,0.00){\vector(1,0){20.00}}
\put(25.00,80.00){\vector(1,0){20.00}}
\put(55.00,0.00){\makebox(0.00,0.00){$\cdots$}}
\put(55.00,80.00){\makebox(0.00,0.00){$\cdots$}}
\put(75.00,90.00){\makebox(0.00,0.00){$d$}}
\put(75.00,10.00){\makebox(0.00,0.00){$d$}}
\put(155.00,10.00){\makebox(0.00,0.00){$d$}}
\put(155.00,90.00){\makebox(0.00,0.00){$d$}}
\put(125.00,40.00){\makebox(0.00,0.00)[l]{$A_{N-1}$}}
\put(205.00,40.00){\makebox(0.00,0.00)[l]{$A_N$}}
\put(115.00,30.00){\vector(0,-1){10.00}}
\put(65.00,80.00){\vector(1,0){20.00}}
\put(65.00,0.00){\vector(1,0){20.00}}
\put(145.00,80.00){\vector(1,0){20.00}}
\put(145.00,0.00){\vector(1,0){20.00}}
\put(115.00,80.00){\makebox(0.00,0.00){$\dr {N-1}N$}}
\put(115.00,0.00){\makebox(0.00,0.00){$\dr {N-1}N$}}
\put(195.00,0.00){\makebox(0.00,0.00){$\dr {N}N$}}
\put(195.00,80.00){\makebox(0.00,0.00){$\dr {N}N$}}
\end{picture}}

\end{center}
As in the proof of Corollary~\ref{rollex}, we apply
Proposition~\ref{3}
to $A_N$. Then for
$$d^{(i)} := d|_{S^{(i)}} \mbox{ for $S^{(i)}:=
\mbox{Span}_{\ci N}(dx^1 \wedge
\cdots \missing{dx^i} \cdots\wedge dx^n) \subset \dr {N-1}N$},
1\leq i\leq
N,$$
we have
\[
A_N\circ d^{(i)} = d\circ\tilde A\circ d^{(i)}+ \chi(A_N)\circ
d^{(i)}.
\]
As before, we
conclude that if the lift $A_{N-1}$ exists, then
$$0=\chi(d\circ A_{N-1}|_{S^{(i)}})=\chi(d\circ \tilde A\circ d^{(i)})+
\chi(\chi(A_N)\circ d^{(i)})=\fr{\pa \chi(A_N)}{\pa x_i}$$
 for each $i$, so  $\chi(A_N)$ is constant.

If $\chi(A_N)$ is constant, then setting $A_{N-1}:= \tilde A d+ 
\chi(A_N)$ and
$A_s := \chi(A_N)$, for $s \leq N-2$, defines a lift
of $A_N$ with the desired property~(\ref{favorit}).%
\qed

\section{More serious models}
\label{ra}

Let $E\to M$ be a smooth vector bundle over a smooth manifold
$M$. We
will, in fact, always suppose that $M = {\bf R}^N$,
with coordinates $\x =
(\rada {x^1}{x^N})$, and that $E$ is the trivial one-dimensional
bundle, $E = \bfr^N \times \bfr \to \bfr^N$,
with only one `vertical' coordinate $u$. More vertical coordinates
present only notational difficulties and all our results directly
generalize to this situation. The case of a general manifold $M$
and a possibly nontrivial bundle $E$ can be studied by standard
 globalization techniques, where our situation will serve as the local
model.

We will consider forms and functions on the infinite jet bundle
$J^\infty E$ over $E$. This jet bundle has
coordinates $ (\x,\u)=(x^i,u_J)$,
where $1\leq i \leq N$, $J$
runs over all multi-indices $(\rada{j_1}{j_N})$,
$\rada{j_1}{j_N}\geq 0$,
and $u_J$ is the coordinate such that
\[
u_J(j^{\infty}(\phi))=\left(\frac{\pa}{\pa \x}\right)^J \phi,
\]
see~(\ref{mx}) for the notation. The order of $J$ is
defined as $|J|=j_1+\cdots+j_N.$

Recall that a {\em local function\/}, $f= f(x^i,u_J)$, 
is by definition the pullback of a smooth function on some $J^k E $, 
and thus depends only on finitely many $u_J$'s. We denote by
$\loc(E)$ the vector space of all local functions. Let
$(\Omega^*(J^\infty E),d)$ be the complex of de~Rham forms on
$J^\infty E$ whose coefficients are local functions. It is
well-known~\cite{anderson:CM92} that the differential on
$\Omega^*(J^\infty E)$ decomposes into a horizontal and a vertical
component  $d=d_H + d_V$, defining the structure of a bicomplex,
the so-called {\em variational bicomplex\/},
\[
\Omega^*(J^\infty E)
= \bigoplus_{k+l=*,\ k,l\geq 0} \Omega^{k,l}(J^\infty E);\
d= d_H+d_V.
\]
Let us denote by $\fr d{dx^i}$,
$1\leq i \leq N$, the total derivative with respect to $x^i$,
\[
\fr{d}{dx^i} := \fr{\pa}{\pa x^i} + \sum_J u_{iJ}\fr{\pa}{\pa u_J},
\]
where
\begin{equation}
\label{Za}
iJ = (\rada {j_1}{j_{i-1}},j_i +1,j_{i+1},\ldots,j_N).
\end{equation}
Given that the $i$-th slot of the multi-index $J$ indicates the number of
$x_i$ derivatives, it would make more sense to
denote one more  $x_i$ derivative on
 $u_J$ by $u_{J+\delta_i}$. Our convention is a (perhaps futile) attempt
to simplify an increasingly complicated system of notation. Let us
remark also that $\fr d{dx^i}$ is usually denoted by $D_i$.  

Then the `horizontal'
differential $d_H : \Omega^{k,*}(J^\infty E)\to
\Omega^{k+1,*}(J^\infty E)$
is given by the formula
\[
d_H\omega =\sum_{1\leq i\leq N} dx^i\wedge\fr{d}{dx^i}\omega .
\]
As in~(\ref{mx}) we denote
\[
\left(\frac d{d\x}\right)^I =
\left(\frac{d}{dx^1}\right)^{i_1} \cdots
\left(\frac{d}{dx^N}\right)^{i_N},
\mbox{ for $I = (\rada {i_1}{i_N})$}.
\]
In order to deal with vertical derivatives, we introduce the
expression
\[
\left(\frac\pa{\pa \u}\right)^\alpha =\prod_J \fr{\pa^{\alpha(J)}}
{\pa u_J^{\alpha(J)}},
\]
where $\alpha$ is a non-negative integer valued function on the
multi-indices $J$ and $\alpha(J)\neq 0$ for only finitely many $J$.

The following definitions are crucial.
\begin{definition}
\label{sa}
The formal  differential degree of
$\alpha$, denoted $\deg_f(\alpha)$, is the
maximal order of multi-index $J$ such that $\alpha(J)\neq 0.$
We say that a linear map $A: \loc(E)\to \loc(E)$ is a LDO
(local differential operator) if it is of the form
\[
\label{sii}
A(f) = \sum_{I,\alpha} p_{I,\alpha}({\x,\u})
\left(\frac{d}{d{\x}}\right)^I
\left(\frac\pa{\pa \u}\right)^\alpha(f),\
f\in \loc(E),
\]
where $p_{I,\alpha}(\x,\u)\in \loc(E)$, and has the property that
\begin{equation}
\label{dfsd}
\mbox{ for each integer $n$, there are only finitely many
$p_{I,\alpha}\not=0$ with  $\deg_f(\alpha)\leq n$}
\end{equation}
\end{definition}
\noindent
Condition~(\ref{dfsd}) guarantees that $A(f)$ is a finite sum
for each $f\in
\loc(E)$. The term `local differential operator' expresses
the  fact that these operators preserve the space of local functions.

It will be useful to introduce the total symbol of a LDO, using
variables $\xi^i, \eta_J$ to represent $\fr d{dx^i}, \fr \pa{\pa u_J}$
respectively. For a LDO $A$, as defined above, we have
\[
\label{symbol}
\sigma(A) = \sum_{I,\alpha} p_{I,\alpha}(\x,\u) (\xi)^I
(\eta)^\alpha
\]

For a LDO $A$ we define the
{\em characteristic\/} to be the {\em differential operator\/} (in
contrast with the case of no vertical variable, where it was
a function)
\begin{equation}
\label{dumka}
\chi(A) := \sum_{I,\alpha}(-1)^I
\left(\left(\frac{d}{d\x}\right)^Ip_{I,\alpha}(\x,\u)\right)
\left(\frac{\pa}{\pa \u}\right)^\alpha \in \ldo.
\end{equation}
Observe that all the horizontal derivatives 
appear only in the coefficients,  
$({d}/{d\x})^Ip_{I,\alpha}(\x,\u),$
so as an operator $\chi(A)$  contains only
vertical derivatives. We will be interested in the lifting problem for
the bottom row of the variational
bicomplex:

\begin{equation}
\el0 \stackrel{d_H}{\longrightarrow}
\el1 \stackrel{d_H}{\longrightarrow} \cdots
\stackrel{d_H}{\longrightarrow}
\el{N-1} \stackrel{d_H}{\longrightarrow}
\el N.
\end{equation}
which is the initial segment of the {\em
Euler-Lagrange complex\/}~\cite{anderson:CM92}.

We extend Definition~\ref{sa} to maps of forms as in
Definition~\ref{bednicky}. Namely, a local
differential operator $A: \el k \to \el l$ is a
linear map whose `matrix coefficients' are LDOs in the sense of
Definition~\ref{sa}.

Then $f \mapsto f \cdot dx^1\land \cdots \land dx^N$
gives an identification $\loc(E)\cong \el N$.
Formula~(\ref{dumka}) thus defines the characteristic also for a
LDO $A: \el N \to \el N$.
The following Proposition is an analog of
Proposition~\ref{3}.

\begin{proposition}
\label{jsem}
For each LDO $A: \el N \to \el N$, there exists a LDO $\tilde A:
\el N \to \el{N-1}$ such that
\[
A = d_H \tilde A + \chi(A).
\]
\end{proposition}

\begin{corollary}
\label{postavil}
A LDO $A_N : \el N\to \el N$ can be lifted into a sequence
$\{A_s : \el s\to \el s\}_{0\leq s\leq N}$ of LDO's
if and only if
\begin{equation}
\label{bub}
\chi(A_N \fr {d}{dx^i})= 0,\
1\leq i\leq N.
\end{equation}
In this case the lift can be chosen in such a way that
\begin{equation}
\label{dnes}
A_s = \chi(A_N),
\mbox { for $0\leq s\leq N-2$.}
\end{equation}
\end{corollary}

\noindent
Here~(\ref{dnes}) means that $A_s$ acts `diagonally' by $A_s(f
(d\x)^\epsilon )= \chi(A_N)(f (d\x)^\epsilon) $.
The corollary will be a consequence of more general statements
of Section~\ref{main_results}.

We end this section with some calculations useful in the
sequel. First, the commutation relation between
$\fr{\pa}{\pa u_J}$
and $\fr{d}{dx^i}$ is given by:
\begin{equation}
\label{snad}
\fr{\pa}{\pa u_J} \fr{d}{dx^i}-\fr{d}{dx^i}\fr{\pa}{\pa u_J}=
\cases{0}{for $j_i =0$, and}{\fr{\pa}{\pa u_K}}{for $J=iK$.}
\end{equation}
It will be convenient to write $\eta_{J/i}$ for $\eta_K$ when  $J=iK$
and, if $i$ does not appear in $J$, then $\eta_{J/i}=0$. 
Relations~(\ref{snad}) can be written very compactly in terms of symbols.
Defining the operator
\begin{equation}\label{theta}
\Theta^i=\sum_J \eta_{J/i}\fr \pa{\pa \eta_J},
\end{equation}
acting as a derivation on symbols. In the symbol calculus, 
commutation relation (\ref{snad}) becomes
\begin{equation}\label{comop}
\eta_J\xi^i-\xi^i\eta_J=\Theta^i(\eta_J).
\end{equation}
Thus for any monomial $\eta^\alpha,$
\[\eta^\alpha\xi^i-\xi^i\eta^\alpha =\Theta^i(\eta ^\alpha).
\]
For reference, we state another commutation relation
in the symbol calculus which  we will need later:
\begin{equation}\label{commutedx}
\xi^ip(\x,\u)=\fr d{dx^i}p(\x,\u)+ p(\x,\u)\xi^i.
\end{equation}

\noindent
Using relations~(\ref{comop})
and~(\ref{commutedx}) we easily deduce that,
\begin{equation}
\label{85}
\sigma(A \fr d{dx^i}- \fr d{dx^i}A) =
(\Theta^i -\fr d{dx^i})\sigma(A),
\mbox { for $A\in \ldo$}.
\end{equation}
The lemma follows immediately from the formula
above.

Define a  `diagonal LDO' $A: \el k \to \el k$ to be 
to be an operator of the form 
$A(f (d \x)^\epsilon)=A(f)(d \x)^\epsilon$, 
where $A$ is a  LDO as above.
\begin{lemma}
A `diagonal LDO' $A$
commutes with the horizontal differential $d_H$ if and only if
\[ (\Theta^i -\fr d{dx^i})\sigma(A)=0,\
\mbox{for all $1\leq i \leq N$}.
\]
\end{lemma}

\noindent
Finally observe that
\[\sigma(\chi(A\fr d{dx^i})) =
(\Theta^i -\fr d{dx^i})\sigma(\chi(A)),\
\mbox{for all $1\leq i \leq N$}.\]
Thus, the assumption (\ref{bub}) implies that
the operators $A_s = \chi(A_N)$ of~(\ref{dnes}) commute with the
differentials. Explicitly:

\begin{proposition}
\label{Benzina}
Suppose $A$ is a LDO as above and that $\chi(A\, \fr d{dx^i})=0$, for
$1\leq i \leq N$. Then the diagonal operator $\chi(A)$ commutes
with
the differential $d_H$.
\end{proposition}
\section{Main results}
\label{main_results}

Let us consider, as in Section~\ref{ra},
the infinite jets on the one-dimensional trivial vector bundle over
$\bfr^N$. First, we need to introduce multilinear differential
operators.

\begin{definition}
\label{hraji}
An {\em $n$-multilinear LDO}
(local differential operator) is an $n$-linear map
$A : {\loc(E)}^{\ot n} \to \loc(E)$ of the form
\begin{equation}
\label{velmi}
A(\rada {f_1}{f_n})\! = \! \sum
p_{I_1,\ldots,I_n,\alpha_1,\ldots,\alpha_n}\!
\left(\left(\frac{d}{d\x}\right)^{I_1}
\left(\frac{\pa}{\pa \u}\right)^{\alpha_1}\!f_1\right)
\!\cdots\!
\left(\left(\frac{d}{d\x}\right)^{I_n}
\left(\frac{\pa }{\pa \u}\right)^{\alpha_n}\!f_n\right)\!,
\end{equation}
where $\Rada f1n\! \in\! \loc(E)$ and
$p_{I_1,\ldots,I_n,\alpha_1,\ldots,\alpha_n} =
p_{I_1,\ldots,I_n,\alpha_1,\ldots,\alpha_n}(\x,\u)\in \loc(E)$
are local functions.
We also require that,
for any $m$, there are only finitely many
multi-indices, $I_1,\dots,I_n;\alpha_1,\ldots,\alpha_n$, such that
\begin{equation}
\label{Pejsek_a_Kocicka}
\sum_{1\leq i \leq n} \deg_f(\alpha_i)\leq m
\mbox{ and }
p_{I_1,\ldots,I_n,\alpha_1,\ldots,\alpha_n}\not=0.
\end{equation}
\end{definition}

\noindent
Condition (\ref{Pejsek_a_Kocicka}) which is the analog of
(\ref{dfsd}) guarantees that the operator
has well defined values on $n$-tuples of local functions.

If we denote by $\fr{d}{d\x_j}$ (resp $\fr{\pa}{\pa \u_j}$)
the total derivative (resp.~the partial derivative)
acting on the $j$-th
function, $1\leq j \leq n$, then we can write the operator
in~(\ref{velmi}) in a more concise form as
\[
A = \sum
p_{I_1,\ldots,I_n,\alpha_1,\ldots,\alpha_n}
\left(\frac{d}{d\x_1}\right)^{I_1}\cdots
\left(\frac{d}{d\x_n}\right)^{I_n}
\left(\frac{\pa}{\pa \u_1}\right)^{\alpha_1}\cdots
\left(\frac{\pa}{\pa \u_n}\right)^{\alpha_n}.
\]
We denote the
vector space of all such $n$-linear local differential operators by
$\ldo(n)$.

\begin{proposition}
\label{vyskop}
The collection $\ldo = \coll{\ldo}$ with the composition maps
\[
\gamma: \ldo(l) \ot \ldo({k_1})\ot \cdots \ot
\ldo({k_1}) \to \ldo(k_1+\cdots+ k_l)
\]
given by $\gamma(A;\rada{A_1}{A_l}): = A(\rada{A_1}{A_l})$
and the
action of the symmetric group  given by $\sigma A(\rada {f_1}{f_n})
:=
A(\rada{f_{\sigma^{-1}(1)}}{f_{\sigma^{-1}(n)}})$,
$\sigma\in \Sigma_n$,
forms an operad.
\end{proposition}

\noindent
{\bf Proof.}
The claim is almost obvious. The only thing which has to be
verified
is that the composition $A(\rada{A_1}{A_l})$ is again a local
differential
operator. But the commutation relation~(\ref{snad}) says how to
move the total derivatives $d/d\x$
over the horizontal derivatives
$\fr{\pa}{\pa u_J}$ to the left, which enables us to write the
composition
$A(\rada{A_1}{A_l})$ in the form~(\ref{velmi}).%
\qed

Define the change of variables,
\[
(\rada
{\x_1}{\x_n})\longmapsto
(\y_1:=\x_1,\y_2:=\x_2-\x_1,\ldots, \y_n:=\x_n-\x_1),
\]
that is, $y^i_1=x^i_1$,
and $y^i_j=x^i_j-x^i_1$ for $2\leq j\leq n$ and $1\leq i
\leq N.$ We have
\begin{equation}
\label{chs}
\fr{d}{d\y_1}
:= \fr{d}{d\x_1} +\cdots + \fr{d}{d\x_n}\quad
\mbox{and}\quad
\fr{d}{d\y_j}
:= \fr{d}{d\x_j},\
\mbox{for $2\leq j\leq n$.}
\end{equation}
Then such an operator, $A$, can be written in the {\em polarized
form\/} (with
different indexing) as
\begin{equation}
\label{pekne}
\sum
q_{I_1,I_2,\ldots,I_n;\alpha_1,\ldots,\alpha_n}
\left(\frac{d}{d\y_1}\right)^{I_1}
\left(\frac{d}{d\y_2}\right)^{I_2}
\cdots
\left(\frac{d}{d\y_n}\right)^{I_n}
\left(\frac{\pa}{\pa \u_1}\right)^{\alpha_1}
\cdots \left(\frac{\pa}{\pa\u_n}\right)^{\alpha_n}.
\end{equation}

\begin{definition}\label{defchar}
For $A \in \ldo(n)$ in the polarized form of~(\ref{pekne}), the
{\em characteristic\/} $\chi(A)\in \ldo(n)$ is defined as
\[
\chi(A) =
 \sum(-1)^{I_1}
\left(\left(\frac d{d\x}\right)^{I_1}
q_{I_1,I_2,\ldots,I_{n};\alpha_1,\ldots,\alpha_n}\right)
\bleft\fr{d}{d\y_2}\bright^{I_2}\!\cdots\!
\bleft\fr{d}{d\y_n} \bright^{I_{n}}
\left(\frac{\pa}{\pa\u_1}\right)^{\alpha_1}
\!\cdots\! \left(\frac{\pa}{\pa\u_n}\right)^{\alpha_n}.
\]
\end{definition}

\noindent
{\bf Remark.} It was necessary to introduce the polarized form 
before defining the characteristic so that  composition with $\fr d{dx^i}$
would introduce only one new  horizontal derivative, $\fr d{dy_1^i}$,
or equivalently, integration by parts on the output of the  multilinear
operator would affect only one tensor component. The relevant formula is 
\begin{eqnarray}
\frac d{dx^i}\circ A&=&\sum
\frac {dq_{I_1,I_2,\ldots,I_{n};\alpha_1,\ldots,\alpha_n}}{dx^i}
\left(\frac{d}{d\y_1}\right)^{I_1}
\left(\frac{d}{d\y_2}\right)^{I_2}
\!\!\ccdots
\left(\frac{d}{d\y_n}\right)^{I_{n}}
\left(\frac{\pa}{\pa \u_1}\right)^{\alpha_1}
\!\!\ccdots 
\left(\frac{\pa}{\pa\u_n}\right)^{\alpha_n}\label{poldef}
\\ 
\nonumber 
&+&\sum q_{I_1,I_2,\ldots,I_{n};\alpha_1,\ldots,\alpha_n}
\left(\frac{d}{d\y_1}\right)^{iI_1}\left(\frac{d}{d\y_2}\right)^{I_2}
\!\!\ccdots 
\left(\frac{d}{d\y_n} \right)^{I_{n}}
\left(\frac{\pa}{\pa \u_1}\right)^{\alpha_1}
\!\!\ccdots \left(\frac{\pa}{\pa\u_n}\right)^{\alpha_n}.
\end{eqnarray}

Let $\ldo^0(n)$ denote the image of $\chi : \ldo(n)\to \ldo(n)$. It
consists of those $A \in \ldo(n)$ whose polarized form does not
contain total derivatives $d/d\y_1$. The endomorphism $\chi$ is a 
projection onto  $\ldo^0(n)$, $\chi^2=\chi$, and
\begin{equation}
\label{Erinmore}
\ldo^0(n)
\cong
\fr{\ldo(n)}{\{A \in \ldo(n);\ \chi(A)=0\}}.
\end{equation}

To simplify the degree conventions, we regrade the horizontal
complex
$\el*$ by introducing
\begin{equation}
\label{regraded}
\Om i := \el{N-i},\ 0\leq i \leq N.
\end{equation}
Thus $(\Om *,d_H)$ is now a chain complex, $\deg(d_H)=-1$.

As we deal with increasingly complicated situations, the
definition of the characteristic becomes more complicated,
but  once
 again, given $A : [(\Om*)^{\ot n}]_0 \to \Om{0}$, there exists
a $\tilde A: [(\Om*)^{\ot n}]_0 \to \Om 1$ such that
\begin{equation}
\label{again}
A=d_H\tilde A + \chi(A).
\end{equation}

Let us introduce the differential graded operad $\dend_* =
\coll{\dend_*}$ of local differential operator endomorphisms 
of the (regraded)
horizontal de~Rham complex $\Om*$.
This means that
$\dend_k(n)$ consists of degree
$k$ graded vector space maps
$f:(\Om*)^{\ot n}\to \Om *$ with `matrix coefficients' from
$\ldo(n)$. An element of $\dend_k(n)$ is thus a sequence $f =
\{f_s\}$, with
\begin{equation}
\label{dymka}
f_s : [(\Om*)^{\ot n}]_s \to \Om{s+k}.
\end{equation}
Observe that $f_s$ may be nonzero only for $\max(-k,0) \leq s \leq
\min(Nn,N-k)$.
The differential $\delta$ on $\dend_*$
is given by the usual formula
\begin{equation}
\label{Sanorin1}
(\delta f)_s := d_H f_s - \sgn{\deg(f)}f_{s-1} d^{\ot n}_H,
\end{equation}
where $d^{\ot n}_H$ is the standard extension of $d_H$ to the tensor product.
Thus $\delta f=0$ if and only if $f$ is a chain map. The
composition
maps and the action of the symmetric group are given as in
Proposition~\ref{vyskop};
the arguments that this indeed defines an operad
structure are the same.

In this context, the lifting problem analogous to the one discussed
in the
previous sections can be formulated as an extending a \ldo,
$f_0 : [(\Om*)^{\ot n}]_0 \to \Om{0}$, to a cocycle in
$\dend_0(n)$.
Observe that
$[{\Om *}^{\ot n}]_0$ consists of elements of the form
\[
\orada{(\omega_1  dx^1\land \cdots \land dx^N)}
{(\omega_n dx^1\land \cdots \land dx^N)},\ \omega_i \in \loc(E),\
1\leq i \leq n,
\]
Since both $\Om 0$ and $[{\om *}^{\ot n}]_0$ are rank one modules over
$\loc(E)$,
$\Om 0\cong\loc(E)$ and  $[{\om *}^n]_0 \cong \loc(E)^{\ot n}$,
so $f_0$ can be interpreted as an
element of $\ldo(n)$. The space $\Om 1$ consists of elements
$\sum\omega_i(dx^1 \wedge \cdots \missing{dx^i} \cdots\wedge
dx^N)$,
and thus is a rank $n$ module, isomorphic to $\loc(E)^{\oplus N}$.

{}From the definition,
\[
[(\Om*)^{\ot n}]_1 = \bigoplus_{1\leq j\leq n}
\Om 0 \ot \cdots \ot \Om1  \ot \cdots \ot \Om 0
\]
($\Om 1$ at the $j$th position).
We derive the identification
\begin{equation}
\label{civka}
[(\Om*)^{\ot n}]_1 = \bigoplus_%
\doubless{1\leq i \leq N}{1\leq j \leq n}
\loc(E)^{\ot n}_{i,j}
\cong (\loc(E)^{\ot n})^{\oplus Nn}.
\end{equation}
Then $\fr d{dx^i_j}$ represents the boundary operator on the
component
of $[(\Om*)^{\ot n}]_1$ corresponding to $\loc(E)^{\ot n}_{i,j}$.

As before, our strategy is to invoke (\ref{again}) to find
$\tilde f_0: [(\Om*)^{\ot n}]_0\to \Om{1}$, to construct $f_1$,
the first
stage of the extension,  using $\tilde f_0$, and
continue from there. 

The identity $\chi(\fr d{dx^i}\circ A)=0$  follows from the
definition of the characteristic. 
Right composition of both sides of  equation  (\ref{again}) for the operator 
 $f_0$ with  $\fr d{dx^i_j}$ gives
\[
f_0 \fr d{dx^i_j}=d_H \tilde f_0 \fr d{dx^i_j}
+ \chi(f_0)\fr d{dx^i_j}.
\]
Equation (\ref{poldef}) immediately implies
\begin{equation}
\label{Alicek}
\chi(d_H A)=0,\
\mbox{ for any \ldo\ $A:[(\Om*)^{\ot n}]_0\to \Om 1$.}
\end{equation}
Moreover, the identity
$\chi(f_0 \fr d{dx^i_j})=\chi(\chi(f_0 )\fr d{dx^i_j})$
shows that a  necessary condition for the existence of a lifting is
\begin{equation}
\label{key}
\chi(f_0 \fr d{dx^i_j})=0,\
1\leq i \leq N,\
1\leq j \leq n.
\end{equation}

Understanding equation (\ref{key}) is the first step towards
a complete description of the $0$-cycles in $\dend_*(n)$.
Basically, we can say that in the multilinear situation, $n\geq 2$,
(\ref{key}) implies that
the terms of formal differential degree zero in the characteristic
determine all the terms of higher formal differential degree.
The precise statement requires some preliminaries.

First we extend the definition of
the symbol to multilinear \ldo, using variables $\xi^i_j$
to represent
the derivatives $\fr d{dx^i_j}$ for $1\leq i\leq N,\,
1\leq j \leq n$ and variables $\eta_J^j$ to represent $\fr
\pa{\pa u_J^j}$
 for $J=(j_1,\ldots, j_N)$ and $1 \leq j \leq n$.
Defining operators analogous to those in (\ref{theta})
\begin{equation}\label{theta2}
\Theta^i_j=\sum \eta_{J/i}^j\fr\pa{\pa \eta_J^j}
\end{equation}
we have the commutation relations
\begin{equation}
\label{comop2}
\eta_J^j  \xi ^i_{k}-\xi^i_{k} \eta_J^j
=\Theta^i_{k}(\eta_J^j)=\delta^j_{k}
\Theta^i_{j}(\eta_J^j).
\end{equation}
Since the symbol determines completely the LDO $A$, we can define
a character $\chi'$, mapping symbols to symbols, with the property
\[ \chi'(\sigma(A))=\sigma (\chi(A)).
\]
The commutation relation (\ref{comop2}) implies
\[ \sigma(A\fr d{dx^i_j})=\sigma(A)\xi^i_j=
\xi^i_j*\sigma(A) +\Theta^i_j\sigma(A),
\]
where
\[
\xi^i_j * \left(
p(\x,\u)(\xi_1)^{I_1} \cdots (\xi_n)^{I_n}
(\eta^1)^{\alpha_1} \cdots (\eta^n)^{\alpha_n}
\right)
:= p(\x,\u)\xi^i_j(\xi_1)^{I_1} \cdots (\xi_n)^{I_n}
(\eta^1)^{\alpha_1} \cdots (\eta^n)^{\alpha_n}.
\]
Since $\fr d{d\y_j}=\fr d{d\x_j}$ for $2\leq j\leq n$ we use
the same
symbol $\xi^i_j$ for $\fr d{dy^i_j}$, but we define a  new symbol
$\zeta^i$ corresponding to the operator $\fr d{dy^i}$.

\begin{proposition}
\label{keyprop}
If A is an $n$-linear LDO, for $n\geq 2$,
then $\chi(A\fr d{dx^i_j})=0$ if and only if the symbol character
$\chi'(\sigma (A))$ satisfies, for $1\leq i \leq N$,
the following system of equations:

\begin{eqnarray}
\label{crux}
\sum_{j=1,\ldots,n}\Theta^i_j (\chi'(\sigma (A))
&=&
\frac d{dx^i}\chi'(\sigma(A)
\\
\label{crux2}
\Theta^i_j(\chi'(\sigma(A)))&=& -\xi^i_j*\chi'(\sigma(A)),\
\mbox{for $j\geq 2$.}
\end{eqnarray}
\end{proposition}
\noindent{\bf Proof.}
First we use the commutation relations to rewrite
\[
\chi'(\sigma(A)\xi^i_j)=
\chi'(\xi^i_j*\sigma(A)+\Theta^i_j(\sigma(A)))
=\chi'(\xi^i_j*\sigma(A)) + \chi'(\Theta^i_j(\sigma(A))).
\]
Then using
$\zeta^i = \prada{\xi^i_1}{\xi^i_n}$
and the identities
\[
\chi'(\zeta^i*\sigma(A))=-\fr d{dx^i}\chi'(\sigma(A))\
\mbox{ and }\
\chi'(\Theta^i_j(\sigma(A)))=\Theta^i_j(\chi'(\sigma(A)),
\]
we deduce
\begin{eqnarray}
\label{chiprime}
\chi'(\sigma(A)\zeta^i)&=& \chi'(\zeta^i*\sigma(A))
+\sum_{j=1,\ldots n}\Theta^i_j(\chi'(\sigma(A)))
\\
\nonumber
&=&-\frac d{dx^i}\chi'(\sigma(A))
+\sum_{j=1,\ldots n}\Theta^i_j(\chi'(\sigma(A))),\,\mbox{and}
\\
\nonumber
\chi'(\sigma(A)\xi^i_j)&=& \chi'(\xi^i_j*\sigma(A))+
\Theta^i_j(\chi'(\sigma(A))),\
\mbox{for $2 \leq j \leq n$.}
\end{eqnarray}
If $\chi(A\fr d{dx^i_j})=0$ for $1\leq j\leq n$, then the left
side of each of these
equations is zero
and rewriting the resulting equations gives the equations
in the statement of the proposition.%
\qed

The crucial fact in understanding equations~(\ref{crux})
and~(\ref{crux2})
is that $\Theta^i_j$ is a
derivation taking $\eta^j_J$ to $\eta^j_{J/i}$,
thus lowering formal differential
degree but preserving the total homogeneity in
all the variables $\eta^j_J$.
There is a simple ordering
on the symbol monomials  $(\eta)^\alpha=
\Pi_J(\eta_J)^{\alpha(J)}$  for a linear
LDO, i.e. $n=1.$ It is defined
by first lexicographically ordering the indices $J$, second,  representing
the exponent $\alpha$ by the sequence of values $\{\alpha(J)\}$
(with finitely many non-zero terms), and
third,  using the lexicographical ordering on these sequences,
reading as in Hebrew from right to left, that is, beginning with the
highest nonzero terms. For example, $\eta^n = (\eta_{(\rada00)}^n)$,
of formal differential degree zero, is minimal among terms
of homogeneity $n$ because it corresponds to the sequence
$(n,0,0,\ldots)$ which is
less than any sequence with a nonzero value beyond the
first term. The above simple order induces a
partial order on the
symbol monomials $\Pi_{j,J}(\eta^j_J)^{\alpha_j(J)}$
with the property that
$\Pi_{j,J}(\eta^j_J)^{\alpha_j(J)}$ is minimal among terms of
the same homogeneity if the formal differential degree of all $\alpha_j$'s
is zero.

%Finally, using the lexicographical order
%on $n$-tuples of sequences, we define a simple order on the symbol
%monomials for the multilinear LDO,
%$(\eta)^\alpha=\Pi_{j,J}(\eta^j_J)^{\alpha_j(J)}.$

Relative to this ordering,
the operators $\Theta^i_j$ all have the effect of lowering
the order.
To simplify notation we introduce the symbols
$\hat I=(I_1,\ldots I_n)$ and
$\hat \alpha=(\alpha_1,\ldots, \alpha_n).$
Thus the coefficient $p_{\hat I,\hat \alpha}(\x,\u)$
appears in the expressions  $\Theta^i_j\chi'(\sigma(A))$
on the left of~(\ref{crux}) and~(\ref{crux2})
multiplied by a monomial $\xi^{\hat I}\eta^{\hat\alpha'}$ where
$\hat\alpha'$ has lower order than $\hat \alpha$.
On the right hand side of~(\ref{crux})
the coefficient of $\xi^{\hat I}\eta^{\hat \alpha'}$
is $\frac d{dx^i}p_{\hat I,\hat\alpha'}$
while on the right hand side of~(\ref{crux2}) the coefficient of
$\xi^{\hat I}\eta^{\hat \alpha'}$ consists of terms involving
$p_{\hat
I',\hat \alpha'}$ for values
$\hat I'$ of lower order relative to the natural
lexicographical order on the
indices $\hat I$. An elementary recursion argument
shows that the solution of equations~(\ref{crux}) and~(\ref{crux2})
is determined uniquely by the coefficients of the terms
$\xi^{\hat I}\eta^{\hat\alpha}$ for
minimal $\alpha$, that is the monomials
$\Pi (\xi_j)^{I_j}(\eta^j)^{m_j}$
with no factors $\eta^j_J$ of positive formal differential degree.

For example, consider a line bundle $E\to {\bf R}$ and
 a bilinear LDO, $A:[(\Om*)^{\ot 2}]_0\to \Om 0$, with symbol
\[
\sigma(A)=\sum
p_{i,j,a,b}(x,u)\zeta^i\xi_2^j\eta_a^1\eta_b^2,
\]
that is, first order and of homogeneity one in the $u$ derivatives.
The symbol character is
\[
\chi(\sigma(A))=
\sum (-1)^i \left(\frac d{dx}\right)^i
p_{i,j,a,b}(x,u)
\xi_2^j\eta_a^1\eta_b^2=:
\sum\chi_{a,b}(x,u,\xi_2)\eta_a^1\eta_b^2.
\]
Equation~(\ref{crux}) becomes
\[
\chi_{a+1,b}+\chi_{a,b+1}= \fr {d\chi_{a,b}}{dx}
\]
and equation~(\ref{crux2}) becomes
\[\chi_{a,b+1}=-\xi_2 *\chi_{a,b}.
\]
Clearly $\chi_{0,0}$ determines all the $\chi_{a,b}$ for
$a,b\geq 0$.

\begin{figure}
\begin{center}
{% Picture saved by xtexcad 2.4
\unitlength=1pt
\begin{picture}(190.00,530.00)(100.00,0.00)
\put(180.00,0.00){\makebox(0.00,0.00){$\vdots$}}
\put(180.00,530.00){\makebox(0.00,0.00){$\vdots$}}
\put(180.00,60.00){\vector(0,-1){50.00}}
\put(180.00,520.00){\vector(0,-1){50.00}}
\put(190.00,230.00){\makebox(0.00,0.00){$a_0$}}
\put(130.00,220.00){\makebox(0.00,0.00)[bl]{$\boxed{3}$}}
\put(130.00,310.00){\makebox(0.00,0.00)[bl]{$\boxed{2}$}}
\put(130.00,410.00){\makebox(0.00,0.00)[bl]{$\boxed{1}$}}
\put(90.00,530.00){\makebox(0.00,0.00){$\vdots$}}
\put(90.00,0.00){\makebox(0.00,0.00){$\vdots$}}
\put(130.00,190.00){\makebox(0.00,0.00){$\B_{-1}$}}
\put(120.00,80.00){\vector(1,0){40.00}}
\put(120.00,180.00){\vector(1,0){40.00}}
\put(120.00,270.00){\vector(1,0){40.00}}
\put(120.00,360.00){\vector(1,0){40.00}}
\put(120.00,450.00){\vector(1,0){40.00}}
\put(180.00,180.00){\makebox(0.00,0.00)[l]%
{\hskip-3mm$\ldo(n)^{\oplus Nn}$}}
\put(180.00,270.00){\makebox(0.00,0.00)[l]{\hskip-3mm$\ldo(n)$}}
\put(180.00,360.00){\makebox(0.00,0.00)[l]%
{\hskip-3mm$\{A \in \ldo(n);\ \chi(A)=0\}$}}
\put(100.00,130.00){\makebox(0.00,0.00){$\delta$}}
\put(100.00,230.00){\makebox(0.00,0.00){$\delta$}}
\put(100.00,320.00){\makebox(0.00,0.00){$\delta$}}
\put(100.00,410.00){\makebox(0.00,0.00){$\delta$}}
\put(100.00,500.00){\makebox(0.00,0.00){$\delta$}}
\put(190.00,320.00){\makebox(0.00,0.00){$a_1$}}
\put(130.00,280.00){\makebox(0.00,0.00){$\B_{0}$}}
\put(130.00,370.00){\makebox(0.00,0.00){$\B_1$}}
\put(90.00,60.00){\vector(0,-1){50.00}}
\put(90.00,520.00){\vector(0,-1){50.00}}
\put(180.00,80.00){\makebox(0.00,0.00){$0$}}
\put(180.00,450.00){\makebox(0.00,0.00){$0$}}
\put(90.00,80.00){\makebox(0.00,0.00){$\dend_{-2}(n)$}}
\put(90.00,180.00){\makebox(0.00,0.00){$\dend_{-1}(n)$}}
\put(90.00,270.00){\makebox(0.00,0.00){$\dend_0(n)$}}
\put(90.00,360.00){\makebox(0.00,0.00){$\dend_1(n)$}}
\put(90.00,450.00){\makebox(0.00,0.00){$\dend_2(n)$}}
\put(180.00,160.00){\vector(0,-1){60.00}}
\put(180.00,430.00){\vector(0,-1){50.00}}
\put(180.00,250.00){\vector(0,-1){50.00}}
\put(180.00,340.00){\vector(0,-1){50.00}}
\put(90.00,160.00){\vector(0,-1){60.00}}
\put(90.00,430.00){\vector(0,-1){50.00}}
\put(90.00,250.00){\vector(0,-1){50.00}}
\put(90.00,340.00){\vector(0,-1){50.00}}
\end{picture}}
\end{center}
\caption{Jacob's un{\bf End}ing ladder.\label{1}}
\end{figure}
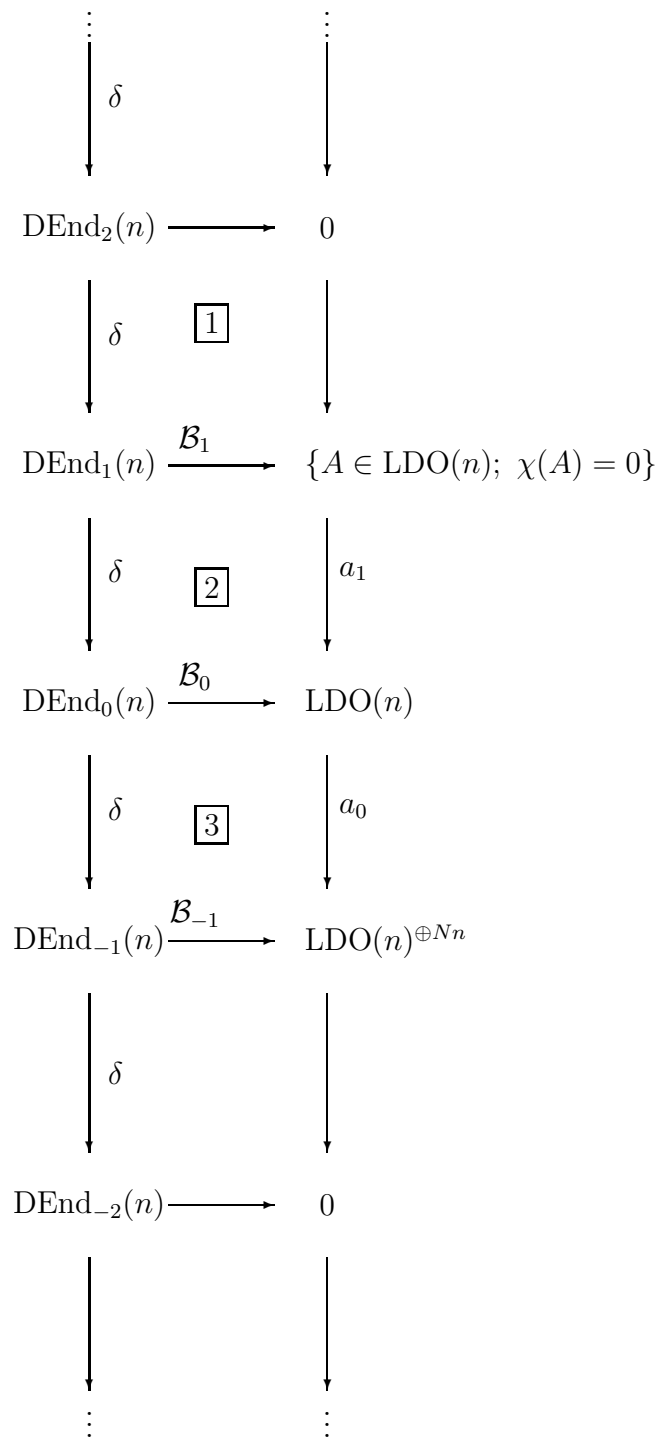

Let us go back to the discussion of the lifting problem.
In Figure~\ref{1} we present a diagram  describing our
situation. It contains  maps $\B_1,\B_0,\B_{-1}$, $a_1$ and
$a_0$ which we now define.

%If we interpret $\ldo$ as a graded operad,
%concentrated in degree zero, then the diagonal map
%$\iota : \ldo \to \dend_*$ of~(\ref{dumka1})
%is a map of graded operads.
For $f=\{f_s\}\in \dend_0(n)$ (the notation of~(\ref{dymka}))
define a projection on the lowest term
 $\B_0 : \dend_0(n)\to \ldo(n)$ by $\B_0(f)
:= f_0$.
For $h=\{h_s\} \in \dend_1(n)$, put $\B_1(h):= d_H
h_0 \in \ldo(n)$. The following lemma is an immediate consequence
of~(\ref{Alicek}).

\begin{lemma}
\label{uhli}
In the situation above, $\chi(d_H h_0) = 0$,
thus $\B_1$ can be interpreted as a map $\B_1: \dend_1(n)\to \{A\in
\ldo; \chi(A)=0\}$.
\end{lemma}

\noindent
The  map  $\B_{-1}:
\dend_{-1}(n)\to\bigoplus_\doubless{1\leq i \leq N}{1\leq j\leq n}
\ldo(n)_{i,j} \cong\ldo(n)^{\oplus Nn} $ is defined as follows. 
Let $g =\{g_s\}\in
\dend_{-1}(n)$, then
\[
g_1 : [(\Om*)^{\ot n}]_1 \to \Om0.
\]
Using the description~(\ref{civka})
of $[(\Om*)^{\ot n}]_1$ as a direct sum  of pieces
$\loc(E)^{\ot n}_{i,j}$,
the $(i,j)$-th component of
$\B_{-1}(g)$ is defined to be the characteristic of
the restriction $g_{1}|_{\loc(E)^{\ot n}_{i,j}}$.

Let  $a_1: \{A\in
\ldo(n);\
\chi(A)=0\}\hookrightarrow \ldo(n)$ be the inclusion and, finally,
 the map $a_0: \ldo(n) \to \ldo(n)^{\oplus Nn}$ is given by
$a_0(A)_{i,j} :=\chi(A\circ \fr {d}{dx_j^i})$.

\begin{lemma}
\label{Praha}
The sequence $\{A\in \ldo(n);\ \chi(A)=0\} \stackrel%
{a_1}{\longrightarrow}
\ldo(n)\stackrel%
{a_0}{\longrightarrow}
\ldo^{\oplus Nn}$
is a differential chain complex.
\end{lemma}

\noindent
{\bf Proof.}
We must show that $a_0a_1=0$, which is the same as to prove that
$\chi(A\circ \frac {d}{dx^i_j})= 0$
whenever $\chi(A)=0$. This
follows immediately from (\ref{chiprime}).%
\qed

\begin{lemma}
All horizontal maps in The un{\bf End}ing ladder (Figure~\ref{1})
are maps of chain complexes.
\end{lemma}

\noindent
{\bf Proof.}
We must show that the diagrams
 \boxed1, \boxed2 and \boxed3 commute.

\noindent
{\em $\boxed1$ commutes.\/}
This means proving that $\B_1\delta(l)=0$, for each $l\in
\dend_2(n)$. But $(\delta l)_0=(d_H l)_0$, thus
$\B_1\delta(l)=0$ follows from $d^2_H=0$.

\noindent
{\em $\boxed2$ commutes.\/}
The equality $a_1\B_1(h)=\B_0\delta(h)$ follows
immediately from definitions.

\noindent
{\em $\boxed3$ commutes.\/}
If $f = \{f_s\}\in \dend_0(n)$
then, by definition, $[a_0\B_0(f)]_{i,j}= \chi(f_0\circ \fr
{d}{dx^i_j})$. On the other hand, $\B_{-1}\delta (f)$ is
computed as
the characteristic of the restriction of $f_0d_H + d_H f_1$
to $\loc(E)_{i,j}$. The second term gives, by Lemma~\ref{uhli},
zero,
while the first term gives $\chi(f_0\circ \fr
{d}{dx_j^i})$, as it should.
\qed

\begin{theorem}
\label{Mixture}
In the diagram of Figure~\ref{1},
\begin{itemize}
\item[(i)]
the complex $(\dend_*(n),\delta)$ is acyclic in positive
dimensions,
$H_{>0}(\dend_*(n),\delta) = 0$,
%\item[(ii)]
%the map $H_0(\iota): \sdo(n)\to H_0(\dend_*(n),\delta)$
%is an isomorphism, and
\item[(ii)]
the map $\B_*$ induces an isomorphism of the
$0$th homology group,
\[
H_0(\dend(n)) \cong \fr{
\{A\in \ldo(n);\ \chi(A \fr{d}{dx^i_j}) = 0,
\mbox{ for all } 1\leq i \leq N,
1\leq j \leq n \}
}
{
\{A\in \ldo(n);\ \chi(A) = 0\}}.
\]
\end{itemize}
Moreover, the map $\B_0$ is an epimorphism of cycles,
$\B_0(Z_0(\dend(n)))=
{\rm Ker} (a_0)$.
\end{theorem}
\noindent{\bf Proof.} See the next section.

\begin{corollary}
\label{Husicka}
A LDO $A: [\Om*^{\ot n}]_0 \to \Om 0$ can be lifted to a sequence
$f_s: [\Om*^{\ot n}]_s\to \Om s$ with $f_0=A$ if and only if
\begin{equation}
\label{divoka}
\chi(A \fr{d}{dx^i_j})=0,
\end{equation}
for each $1\leq i \leq N$, $1\leq j\leq n$.
In the case $\chi(A)=0$, the lift can be chosen in such a way that
\begin{equation}\label{kapka}
f_s = 0,\ \mbox{for $s\geq 2$.}\end{equation}
\end{corollary}

\noindent
{\bf Proof.}
The first part of the corollary claims the existence of an $f\in
\dend_0(n)$, $\delta f=0$, with $\B_0(f)=A$.
But~(\ref{divoka}) means
that $a_0(A)=0$ and the existence of $f$ follows from the fact that
$\B_0$ is an epimorphism of cycles.

Let us prove the second part of the corollary. Suppose that
$\chi(A)=0$ and let $\overline f=\{\overline f_s\}$ be a lift
of $A$. Since $A\in \mbox{\rm
Im}(a_1)$, there exists $h\in \dend_1(n)$ such that $\overline f =
\delta h$. This means that $A = \overline f_0 = d_H
h_0$. One immediately sees that $f = \{f_s\}$ with $f_0 = A$, $f_1
:= h_0d^{\ot n}_H$ and $f_s=0$ for $s\geq 2$ is a lift of $A$.\qed

%The proof of the general case follows by applying the above on
%$\overline A := A -\chi(A)$, then $\{\chi(A)+ f_s\}$ is a lift of
%$A$ having the desired property~(\ref{kapka}).\qed

\section{Proofs}
\label{proofs}

This section is devoted to the proof of Theorem~\ref{Mixture}. The
basic tool will be the following
de~Rham complex with operator coefficients.

\begin{definition}
The {\em operator complex\/} $O^*(n)= (O^*(n),d)$
is the complex of de~Rham forms on $J^\infty E$ with
coefficients in $\ldo(n)$. The differential $d$ is given by
\begin{equation}
\label{letadylko}
d(A (d\x)^\epsilon) = \sum_{1\leq i\leq N}
\left(\frac d{dx^i}A\right) dx^i\wedge (d\x)^\eps.
\end{equation}
\end{definition}

\noindent
Let $J : O^N(n)\to \ldo(n)$ be the map $J(A{dx^1}
\wedge\cdots\wedge{dx^N})
:= A$.

\begin{theorem}
\label{zapalky}
The complex $(O^*(n),d)$ is acyclic in degrees $< N$, that is
$H^{<N}(O^*(n),d)=0$, while the map $J$ induces an isomorphism
\[
H^N(O^*(n),d) \cong \fr{\ldo(n)}{\{A \in \ldo(n);\ \chi(A)=0\}}.
\]
\end{theorem}

Notice the following rather surprising fact: the operator complex
$(O^*(n),d)$ is acyclic in degree $0$, though the `ordinary'
horizontal
de~Rham
complex $(\el *,d_H)$ is not, $H^0(\el *,d_H)= {\bf R}$!
This follows from (and implies)
the following stunning property of local
differential operators:
\[
\mbox{if $A\in \ldo(n)$ and $\fr d{dx^i} \circ A
= 0$, $1\leq i \leq N$, then
$A=0$.}
\]
This will not be true if we remove the `convergence
property'~(\ref{Pejsek_a_Kocicka}). As an example, take the
operator
\[
A:= 1-x \frac d{dx} +\frac 12 x^2 \frac{d^2}{dx^2} - \frac 16
x^3 \frac
{d^3}{dx^3} + \cdots
\]
in one space and no vertical variables. It clearly satisfies $\fr
d{dx} \circ A = 0$. Observe that, for a polynomial $f$, $A(f)=
f(0)$.

\begin{figure}
\begin{center}
{% Picture saved by xtexcad 2.4
\unitlength=1.500000pt
\begin{picture}(200.00,120.00)(0.00,0.00)
\put(20.00,80.00){\makebox(0.00,0.00){$\bullet$}}
\put(40.00,60.00){\makebox(0.00,0.00){$\bullet$}}
\put(0.00,100.00){\makebox(0.00,0.00){$\cdots$}}
\put(0.00,20.00){\makebox(0.00,0.00){$\cdots$}}
\put(0.00,40.00){\makebox(0.00,0.00){$\cdots$}}
\put(0.00,60.00){\makebox(0.00,0.00){$\cdots$}}
\put(0.00,80.00){\makebox(0.00,0.00){$\cdots$}}
\put(180.00,120.00){\makebox(0.00,0.00){$\vdots$}}
\put(160.00,120.00){\makebox(0.00,0.00){$\vdots$}}
\put(140.00,120.00){\makebox(0.00,0.00){$\vdots$}}
\put(120.00,120.00){\makebox(0.00,0.00){$\vdots$}}
\put(100.00,120.00){\makebox(0.00,0.00){$\vdots$}}
\put(80.00,120.00){\makebox(0.00,0.00){$\vdots$}}
\put(60.00,120.00){\makebox(0.00,0.00){$\vdots$}}
\put(40.00,120.00){\makebox(0.00,0.00){$\vdots$}}
\put(20.00,120.00){\makebox(0.00,0.00){$\vdots$}}
\put(70.00,30.00){\makebox(0.00,0.00)[l]{$d_2$}}
\put(90.00,30.00){\makebox(0.00,0.00)[t]{$d_1$}}
\put(80.00,0.00){\makebox(0.00,0.00)[b]{$p=0$}}
\put(200.00,20.00){\makebox(0.00,0.00)[l]{$q=0$}}
\put(140.00,20.00){\makebox(0.00,0.00){$\bullet$}}
\put(120.00,20.00){\makebox(0.00,0.00){$\bullet$}}
\put(100.00,20.00){\makebox(0.00,0.00){$\bullet$}}
\put(80.00,20.00){\makebox(0.00,0.00){$\bullet$}}
\put(120.00,40.00){\makebox(0.00,0.00){$\bullet$}}
\put(100.00,40.00){\makebox(0.00,0.00){$\bullet$}}
\put(80.00,40.00){\makebox(0.00,0.00){$\bullet$}}
\put(60.00,40.00){\makebox(0.00,0.00){$\bullet$}}
\put(100.00,60.00){\makebox(0.00,0.00){$\bullet$}}
\put(80.00,60.00){\makebox(0.00,0.00){$\bullet$}}
\put(60.00,60.00){\makebox(0.00,0.00){$\bullet$}}
\put(80.00,80.00){\makebox(0.00,0.00){$\bullet$}}
\put(60.00,80.00){\makebox(0.00,0.00){$\bullet$}}
\put(40.00,80.00){\makebox(0.00,0.00){$\bullet$}}
\put(60.00,100.00){\makebox(0.00,0.00){$\bullet$}}
\put(40.00,100.00){\makebox(0.00,0.00){$\bullet$}}
\put(20.00,100.00){\makebox(0.00,0.00){$\bullet$}}
\put(180.00,10.00){\line(0,1){100.00}}
\put(160.00,10.00){\line(0,1){100.00}}
\put(140.00,10.00){\line(0,1){100.00}}
\put(120.00,10.00){\line(0,1){100.00}}
\put(100.00,10.00){\line(0,1){100.00}}
\put(80.00,10.00){\line(0,1){100.00}}
\put(60.00,10.00){\line(0,1){100.00}}
\put(40.00,10.00){\line(0,1){100.00}}
\put(20.00,10.00){\line(0,1){100.00}}
\put(10.00,100.00){\line(1,0){180.00}}
\put(10.00,80.00){\line(1,0){180.00}}
\put(10.00,60.00){\line(1,0){180.00}}
\put(10.00,40.00){\line(1,0){180.00}}
\put(10.00,20.00){\line(1,0){180.00}}
\thicklines
\put(80.00,20.00){\vector(0,1){20.00}}
\put(80.00,20.00){\vector(1,0){20.00}}
\end{picture}}
\end{center}
\caption{The shape of the bicomplex ${\bf E} =
(\EE**,d=d_1+d_2)$ for $N=3$.
Solid dots indicate nontrivial entries.\label{shape}}
\end{figure}
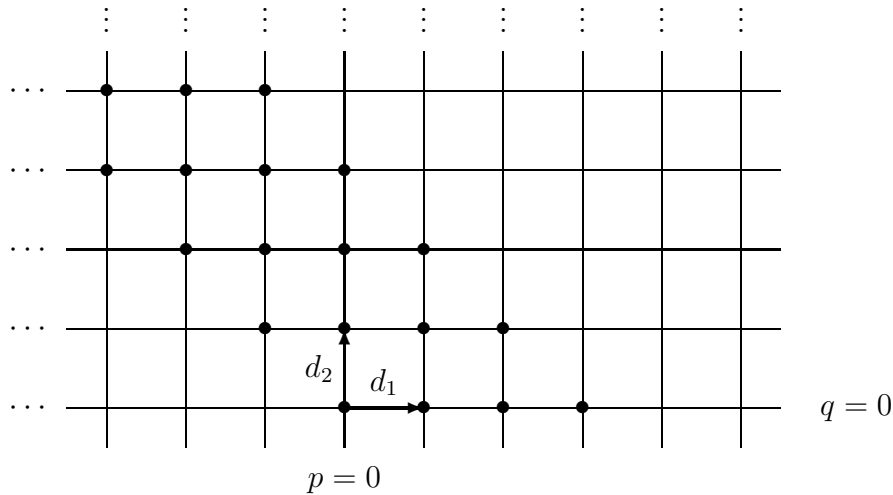

\noindent
{\bf Proof of Theorem~\ref{zapalky}.}
Let us look more closely at the structure of the differential
in the
complex $(O^*(n),d)$.
If $A (d \x)^\eps \in O(n)$, where $A\in \ldo(n)$ is
as in~(\ref{pekne}), then 
\[
d(A(d \x)^\eps ) =(d_1+d_2)(A(d \x)^\eps )
\]
with
\[
d_1(A(d \x)^\eps ) =
\sum
\frac{ dq_{I_1,I_2,\ldots,I_{n};\alpha_1,\ldots,\alpha_n}}{dx^i}
\left(\frac{d}{d\y_1}\right)^I
\left(\frac{d}{d\y_2}\right)^{I_2}
\!\!\!\ccdots
\left(\frac{d}{d\y_n}\right)^{I_{n}}
\left(\frac{\pa}{\pa \u_1}\right)^{\alpha_1}
\!\!\!\ccdots
\left(\frac{\pa}{\pa\u_n}\right)^{\alpha_n}
dx^i\wedge (d \x)^\eps
\]
and
\[
d_2(A(d \x)^\eps ) =
\sum
q_{I_1,I_2,\ldots,I_{n};\alpha_1,\ldots,\alpha_n}
\left(\frac{d}{d\y_1}\right)^{iI}
\left(\frac{d}{d\y_2}\right)^{I_2}
\!\!\!\ccdots
\left(\frac{d}{d\y_n} \right)^{I_{n}}
\left(\frac{\pa}{\pa \u_1}\right)^{\alpha_1}
\!\!\!\ccdots
\left(\frac{\pa}{\pa\u_n}\right)^{\alpha_n}
dx^i\wedge (d \x)^\eps.
\]
Let us decompose $O^*(n) = \bigoplus_{p+q=*}O^{p,q}(n)$, where
$O^{p,q}$
consists of $A(d \x)^\eps  \in O^{p+q}(n)$ such that $A\in \ldo(n)$
in the
polarized form of~(\ref{pekne})
contains exactly $q$ instances of total derivatives $\fr
d{dy}$. Denote for simplicity $\EE pq := O^{p,q}(n)$. Then
\[
d_1 : \EE pq \to \EE{p+1}q
\mbox{ and }
d_2: \EE pq \to \EE p{q+1}.
\]
We are going to study the properties of the bicomplex ${\bf E} :=
(\EE**,d=d_1+d_2)$. Observe that ${\bf E}$
is not contained in the first
quadrant, but
\[
\EE pq \not= 0
\mbox{ for $0\leq p+q \leq N$, $q\geq 0$.}
\]
The shape of the bicomplex is indicated in Figure~\ref{shape}.
For any fixed $q,m\geq 0$, let $V^m_q$ be, , the
${\bf R}$-vector space with the
basis
\[
\{
e^{I,\Rada I2n;\Rada \alpha1n};\
|I|=q \mbox{ and } \sum_{1\leq i \leq n}\deg_f(\alpha_i) \leq m
\}
\]
and let $V_q := \displaystyle{ \lim_{\longleftarrow}}
V^m_q$, the inverse
limit of the system of natural projections $V^m_q \to V^n_q$,
$m\geq n$.
Then the horizontal complex
$(\EE *q,d_1)$ is
the tensor product of $V_q$ and
the horizontal de Rham complex $(\el *,d_H)$. It follows from the
description of the cohomology of this
complex~\cite{anderson:CM92} that
\begin{equation}
\label{Shackelton}
H^s(\EE *q,d_1) =
\tricases {V_q}%
{for $s=-q,$}
{0}{for $-q<s< N-q$, and}%
{V_q \ot H^N(\el *,d_H)}{$s=N-q$.}
\end{equation}

One is tempted to calculate the
homology of the bicomplex ${\bf E}$
using an obvious spectral sequence, but, since the bicomplex
is not only in the first quadrant, one must be very
careful about the convergence.
The calculation just given of the first stage
of the spectral sequence of a double complex
filtered by rows, the complex
with operator $d_1$, gives an $E_1$
term with nonzero entries along both
counterdiagonals $p+q=0$ and $p+q=N$.
It would seem from this that the cohomology of
the double complex is far from
trivial in dimension zero. The reason that this is not the case is
that the bicomplex is not a first quadrant one and using
the filtration by columns would lead us to
constructing a cocycle
which would be an infinite sum with increasing $q$,
contradicting the finiteness condition (\ref{Pejsek_a_Kocicka}).

On the other hand, consider the spectral sequence for the
filtration
by columns.
Let $W^m$ be, for $m\geq 0$, the $\loc(E)$-vector space with
the basis
\[
\{f^{\Rada I2n; \Rada \alpha 1n};\
\sum_{1\leq i \leq n} \deg_f(\alpha_i) \leq m\}
\]
and let $W := \displaystyle{\lim_{\longleftarrow}W^m}$, of
course, $W
\cong \ldo^0(n)$. Then the $p$-th column of the
first stage of this spectral sequence,
with differential given by $d_2$, is isomorphic to $W$ tensored
with
the Koszul complex
\[
K_p:=\cdots
\longrightarrow
S^q(\zeta)\ot \wedge^{p+q} (d\x)\stackrel{d_K}{\longrightarrow}
S^{q+1}(\zeta)\ot \wedge^{p+q+1} (d\x)\longrightarrow\cdots\]
\[d_K(f\ot (d\x)^\eps)=\sum_{i=1,\ldots,N}
\zeta^i f\ot dx^i\wedge (d\x)^\eps
\]
where $0\leq q \leq N-p$, $S^q(\zeta)$ is
the $\bf R$-vector space of homogeneous
polynomials of degree $q$ in the variables
$\zeta^1,\ldots,\zeta^N$, and
$\bigwedge^{p+q}(d\x)$ is the degree $p+q$ component of
the $\bf R$-exterior algebra on $dx^1,\ldots, dx^N$.  It is easy
to see,
for example using the contracting homotopy which is
defined on terms
$f\ot (d\x)^{\epsilon} \in S^q(\zeta)\ot \land^{p+q}(d\x)$ by
\[
f\ot
(d\x)^\eps\longmapsto
\fr1{N-p}\sum_{1\leq i \leq N}
\fr {\pa f}{\pa \zeta^i}\ot
\iota\left(\frac \pa {\pa x^i}\right)(d\x)^\eps,
\]
that the cohomology of the  Koszul complex is trivial for $p< N$
and one dimensional
with basis $dx^1\wedge \cdots \wedge dx^N$ for $p=N$. Thus
\begin{eqnarray*}
H^q(\EE p*,d_2)&=&0, \mbox{ for $0< q$, or $0=q$ and $p< N$, and}
\\
H^0(\EE N*,d_2)&=& W \cong \ldo^0(n).
\end{eqnarray*}
In this case, since $E_1$ has only one nonzero term, we can
conclude
that
$H^N(O^*(n),d) \cong {\ldo^0(n)}$,
but we shall give a direct proof
anyway.

Suppose that $c = c_{k,0}+ c_{k-1,1} + \cdots c_{k-u,u}$ is
a degree
$k$ cycle, $0< k \leq N$, $c_{k-i,i}\in \EE{k-i}i$, $i\geq 0$.
By a degree argument,
$d_2(c_{k-u,u})= 0$. Let $u\geq 1$. By the acyclicity of
$(\EE p*,d_2)$, there exists an $\alpha \in \EE{k-u}{u-1}$
such that
$d_2(\alpha) = c_{k-u,u}$. We then replace $c$
by $c - (d_1+d_2)(\alpha)$,
which is in the same homology class, but which has no component in
$\EE{k-u}u$. Repeating this process as many
times as
necessary we conclude that we could in fact assume that $u=0$,
or $c
\in \EE k0$. This means, since we assumed $c$ to be a cycle, that
$d_1 c=0$ and $d_2 c=0$.
Observe that this reduction works for $k=0$ as
well; we may immediately conclude that $c_{-u,u}=0$ since
$d_2$ is a monomorphism on $\EE{-u}{u}$.

Now, if $k < N$, we immediately conclude
that $c=0$, because $d_2$ is a
monomorphism on $\EE k0$. This proves the acyclicity of $O^*(n)$ in
degrees $< N$.

If $k=N$, then $c = A dx^1\land \cdots\land dx^N$,
where $A \in \ldo(n)$ contains
no $\fr d{dx^i}$; we denoted the set of all such \ldo's by
$\ldo^0(n)$. Such $c$ cannot be a nontrivial boundary. To see this,
let $c = (d_1+d_2)b$, with $b = b_{k-1,0}+ b_{k-2,1}+ \cdots+
b_{k-v-1,v}$, $b_{k-i-1,i} \in E^{k-i-1,i}$.
The `leading' term $b_{k-v-1,v}$ must be a
$d_2$-cycle hence a $d_2$-boundary and we may go `down
the staircase' as above and assume $b = b_{k-1,0}$. Then $d_1 (b)
= c$ and
$d_2(b) = 0$, which implies $b = 0$,
since $d_2$ is a monomorphism on
$E^{k-1,0}$.
We proved
\[
H^N(O^*(n),d) \cong {\ldo^0(n)},
\]
which, together with~(\ref{Erinmore}), finishes the proof.\qed

\noindent
{\bf Proof of Theorem~\ref{Mixture}.}
For simplicity, we will explicitly specify the range of summations
only when it will not be obvious.
Suppose $f\in \dend_k(n)$. This means that $f$ is a sequence
of maps
\[
f_s : [(\om*)^{\ot n}]_s \to \om{s+k}.
\]
The space $[(\om*)^{\ot n}]^s$ is spanned by elements of the form
\[
\orada{(h_1 (d\x_1)^{\eps_1})}{(h_n  (d\x_n)^{\eps_n})},\
h_j \in \ci N,\
1\leq j\leq n,
\]
where the subscript indicates to which copy of ${\bf R}^N$
the corresponding object
applies, and $\sum_{1\leq j \leq n} |\eps_j| = nN-s$
(remember
the regrading~(\ref{regraded})). With this notation,
$f_s$ acts by
\[
f_s(\orada{(h_1(d\x_1)^{\eps_1})}
{(h_n(d\x_n)^{\eps_n})})=
\sum_\eps A_{s,\eps}^{\Rada \eps1n}(\Rada h1n)(d\x)^\eps,
\]
with some $A_{s,\eps}^{\Rada \eps1n}\in \ldo(n)$.
In other words, $f = \{f_s\}$ is represented by the system
\begin{equation}
\label{sesivacka}
\{F^{\Rada \eps1n}_s \in O^{N-s-k}(n)\},\
F^{\Rada \eps1n}_s := \sum_\eps A_{s,\eps}^{\Rada \eps1n}(d
\x)^\eps,
\end{equation}
where $\sum_j |\eps_j| = nN-s$ and
$\max(-k,0) \leq s\leq \min(nN,N-k)$. The last inequality
simplifies,
for $k\geq 0$, to $0\leq s \leq N-k$.

Let us try to understand how the differential $\delta$
in $\dend_*(n)$, defined by~(\ref{Sanorin1}), works.
We have
\begin{eqnarray*}
\lefteqn{\hskip -1cm
(\delta f)_s(\orada{(h_1(d\x_1)^{\eps_1})}{(h_n
(d\x_n)^{\eps_n})})
= \sum_{i,\eps}
\frac d{dx^i}A_{s,\eps}^{\Rada \eps1n}(\Rada h1n)dx^i\wedge
(d\x)^\eps}
\\
&&-  \sum_{i,j,\delta}
\znamenko{k+\prada{|\eps_1|}{|\eps_{j-1|}}}
A_{s,\delta}^{\eps_1,\ldots,i\eps_j,\ldots,\eps_n}
(h_1,\ldots,\frac d{dx^i}h_j,\ldots, h_n)(d\x)^\delta .
\end{eqnarray*}
In the second term, $i\eps_j$ has the obvious meaning similar
to that
of $iJ$, see~(\ref{Za}).
The above formula
can be written in terms of the expressions~(\ref{sesivacka}) as
\begin{equation}
\label{f2}
\delta \{F^{\Rada \eps1n}_s\}
= \{ d F^{\Rada \eps1n}_s -\sum_{i,j}
\znamenko{k+\prada{|\eps_1|}{|\eps_{j-1|}}}
F^{\eps_1,\ldots,i\eps_j,\ldots,\eps_n}_{s-1}\circ
\frac {d}{dx^i_j}
\}
\end{equation}
where the differential $d$ in the first term of the right hand
side is
the differential of the operator complex~(\ref{letadylko}).
Let us prove that {\em $\dend_*(n)$ is acyclic in positive
dimensions.\/}
If $\{F^{\Rada \eps1n}_s\}\in \dend_k(n)$ is a cycle, $k >
0$, then, by~(\ref{f2}),
\begin{equation}
\label{dlachatko}
0= d F^{\Rada \eps1n}_s -\sum_{i,j}
\znamenko{k+\prada{|\eps_1|}{|\eps_{j-1|}}}
F^{\eps_1,\ldots,i\eps_j,\ldots,\eps_n}_{s-1}\circ
\frac {d}{dx^i_j},
\end{equation}
for all $s$, $\Rada {\eps}1n$.
We are looking for $\{H^{\Rada \eps1n}_s \} \in \dend_{k+1}(n)$
such that
\begin{equation}
\label{f1}
F^{\Rada \eps1n}_s = d H^{\Rada \eps1n}_s  +\sum_{i,j}
\znamenko{k+\prada{|\eps_1|}{|\eps_{j-1|}}}
H^{\eps_1,\ldots,i\eps_j,\ldots,\eps_n}_{s-1}\circ \frac
{d}{dx^i_j}.
\end{equation}
Let us solve this equation inductively. For $s=0$ it reduces to
\begin{equation}
\label{Ja}
F^{\Rada \eps1n}_0 = d H^{\Rada \eps1n}_0 ,
\end{equation}
which must be solved in $O^{N-k}(n)$. Equation~(\ref{dlachatko})
with
$s=0$ says
that $d F^{\Rada \eps1n}_0 =0$ and the existence of $H^{\Rada
\eps1n}_0$  follows from the acyclicity of
the operator complex (Theorem~\ref{zapalky}).

Suppose we have solved~(\ref{f1}) for all $s < r$ and try to
solve it
for $s=r$. By the acyclicity of the operator complex in dimension
$N-r-k$,
it is enough to verify that
\begin{equation}
\label{kos1}
F^{\Rada \eps1n}_r  -\sum_{i,j}
\znamenko{k+\prada{|\eps_1|}{|\eps_{j-1|}}}
H^{\eps_1,\ldots,i\eps_j,\ldots,\eps_n}_{r-1}\circ \frac
{d}{dx^i_j}
\end{equation}
is closed in $O^{N-r-k}(n)$. By~(\ref{dlachatko}),
$d F^{\Rada \eps1n}_r$
equals
\[
\sum_{i,j}\znamenko{k+\rada{|\eps_1|}{|\eps_{j-1}|}}
F^{\eps_1,\ldots,i\eps_j,\ldots,\eps_n}_{r-1}\circ
\frac {d}{dx^i_j},
\]
while the differential of the second term of~(\ref{kos1}) equals,
by the inductive assumption, to
\begin{equation}
\label{Kveta}
-\!\sum_{i,j}\!
\znamenko{k+\prada{|\eps_1|}{|\eps_{j-1|}}}
F^{\eps_1,\ldots,i\eps_j,\ldots,\eps_n}_{r-1}\!\circ\!
\frac {d}{dx^i_j}+\hskip -3mm
\sum_\doubless{i,j,k,l}{j\not=l}%
(-1)^{\sigma(j,l)}\!\cdot\!
H_{r-2}^{\eps_1,\ldots,i\eps_j,\ldots,k\eps_l,\ldots,\eps_n}
\!\circ\! \frac {d}{dx^i_j}\frac{d}{dx^k_l}
\end{equation}
where
\[
\sigma(j,l) :=
\cases {1+\prada{|\eps_j|}{|\eps_{l-1}|}}{for $j< l$, and}%
{\prada{|\eps_l|}{|\eps_{j-1}|}}{for $j> l$.}
\]
It is immediate to conclude that the second term of~(\ref{Kveta})
is
zero. Summing up the above informations we
see that the form in~(\ref{kos1}) is indeed closed,
and the induction may go on.

Let us assume $f = \{F^{\Rada\eps1n}_s\}\in \dend_0(n)$. This means
that $F^{\Rada\eps1n}_s
\in O^{N-s}(n)$, $0\leq s\leq N$, and $\sum_j |\eps_j|
= nN-s$. For $s=0$ this may happen only if $\eps_j = (\rada11)$
for all
$1\leq j\leq N$,
and the system $\{F^{\Rada\eps1n}_0\}$ boils down to one element
$F_0 \in O^N(n)$. Let us write
$F_0 = f_0 dx^1 \land \ldots \land dx^N$ with some
$f_0 \in \ldo(n)$. By definition, $\B_0(f) = f_0$. As before,
try to
solve~(\ref{f1}) inductively. For $s=0$ it reduces to
\[
F_0 = dH_0,\ H_0 \in O^{N-1}(n),
\]
which can be, by Theorem~\ref{zapalky}, solved if and only if
$\chi(f_0)=0$, which is the same as
$\B_0(f) \in {\rm Im}(a_1)$. If this is the case,
the induction may go on, by the acyclicity of the operator complex.
We proved that {\em $H_0(\B_*)$ is a monomorphism.\/}

Let us prove that {\em$\B_0(Z_0(\ldo(n)))=
{\rm Ker} (a_0)$.\/} Suppose $f_0 \in \ldo(n)$ with $a_0(f_0)=0$. We
are
looking for a cycle $f = \{F_s^{\Rada \eps1n}\}\in \dend_0(n)$ such
that $F_0 = f_0 dx^1 \land \ldots \land dx^N$. We construct such
a cycle by inductively solving~(\ref{dlachatko}):
\begin{equation}
\label{pero}
dF_s^{\Rada \eps1n} =
\sum_{i,j}\znamenko{|\eps_1|+\cdots+|\eps_{j-1}|}
F_{s-1}^{\eps_1,\ldots,i\eps_j,\ldots,\eps_n}\circ
\frac{d}{dx^i_j},
\end{equation}
in $O^{N-s+1}(n)$, for $s\geq 1$.
We already observed that
$F_0^{\Rada \eps1n}$ reduces to one element,
$F_0$. Thus~(\ref{pero}) can be, for $s=1$, written as
\begin{equation}
\label{Kralicek}
dH^{i,j} = (-1)^{N(j-1)}\cdot
F_0 \circ \frac {d}{dx^i_j},\ 1\leq i \leq N,\
1\leq j \leq n,
\end{equation}
where $H^{ij} : =  F_1^{\Rada \eps1n}$ with
\[
\eps_k :=
\cases{(\rada 11)}{ for $k\not= j$, and}{(\rada11,0,\rada11)
\mbox{ /$0$ at the $i$-th position/}}{for $k=j$.}
\]
By Theorem~\ref{zapalky}, equation~(\ref{Kralicek}) can be
solved if
and only if
$\chi( F_0\circ
\fr{d}{dx^i_j}) = 0$, $1\leq i \leq N$, $1\leq j \leq n$,
which is the same as to say that $f_0 \in {\em Ker}(a_0)$.
The inductive construction then may go on by the acyclicity of the
operator complex.\qed

\section{Applications}
\label{applications}
This paper originated in our attempts to understand the
paper~\cite{barnich-fulp-lada-stasheff:preprint}
with the problem of defining a lifting of a Poisson structure
on the
algebra of functionals to a strong homotopy Lie structure on
the horizontal
complex. We review the situation discussed in that paper.
First let us recall a basic definition and
lemma from~\cite{barnich-fulp-lada-stasheff:preprint}.

\begin{definition}
\label{localfunctional}
A {\em local functional}
\begin{eqnarray}
{\cal L}[\phi]=\int L(\x,\phi(\x))dx^1\wedge\cdots\wedge dx^N =
\int (j^\infty\phi)^*  L(\x,\u)dx^1\wedge\cdots\wedge dx^N
\end{eqnarray}
is the integral of the pull-back of an element of 
$L(\x,\u)dx^1\wedge\cdots\wedge dx^N \in \Omega^{N,0}(J^\infty E)$
by the  section $j^\infty \phi$ of $J^\infty E$ corresponding to the section
$\phi$ of $E$, where we assume that $L(\x,0)=0$. 
The integral will always be well-defined since 
$(j^\infty \phi)^*L(\x,\u)dx^1\wedge\cdots\wedge dx^N$ is a compactly
supported smooth $N$-form on ${\bf R}^N$ for $\phi$ a
section with compact support, 
since the coefficient vanishes outside the support
of $\phi.$

The {\em space of local functionals} ${\cal F}$ is the vector
space of equivalence classes of local functionals, where two local
functionals are  equivalent if they agree for all sections of compact support.
\end{definition}

\noindent
The correspondence
between the functional and the corresponding
element of $\Omega^{N,0}(J^\infty E)$ is not
one to one, as we see from the following lemma.

\begin{lemma}
\label{l4}
The vector space of local functionals ${\cal F}$ is isomorphic
to the
cohomology group
\[
H^N(\om{*,0},d_H)=\om{N,0}/d_H(\om{N-1,0}).
%=\om{N,0}/Ker E,
\]
\end{lemma}
The isomorphism of the lemma is induced by the correspondence 
\[
\om{N,0}
\ni f(\x,\u) dx^1 \land \cdots \land dx^N 
\longleftrightarrow {\cal L}
\in {\cal F}, 
\]
where ${\cal L}$ is the functional corresponding to the
form $ L(\x,\u) dx^1 \land \cdots \land dx^N$, with $L(\x,\u) :=
f(\x,\u)- f(\x,0)$.

\noindent
We cite from the introduction
to~\cite{barnich-fulp-lada-stasheff:preprint}:

{\em
``The approach to Poisson brackets in this context, pioneered by
Gel'fand, Dickey and Dorfman
(see \cite{gel1,gel2,gel3,gel4,gel5,Olver,dickey} for
reviews),  is to consider the Poisson brackets for
local functionals as being induced by brackets for local functions,
which are not necessarily strictly Poisson. We will analyze here in
detail the structure of the brackets for local functions
corresponding
to the Poisson brackets for local functionals. More precisely,
we will
show that these brackets will imply higher order brackets combining
into a strong homotopy Lie algebra.''
}

A suitable bracket for local functions,
which is not necessarily strictly Poisson,
is given by $\tilde l_2\in O^N(2)$
considered as a map from $\om{N,0}\otimes
\om{N,0}$ to $\om{N,0}$ such that:

\begin{itemize}
\item[(i)]
$ \tilde l_2(\alpha,d_H\beta)\in d_H\om{N-1,0}$ for
all $\alpha\in\om{N,0}$ and $\beta\in \om{N-1,0}$,
\item[(ii)]
$\tilde l_2(\alpha, \beta)+
\tilde l_2(\beta,\alpha)\in d_H(\om{N-1,0})$, and
\item[(iii)]
\begin{eqnarray*}
\lefteqn{\hskip-1cm
\sum_{\sigma\in \unsh(2,1)}(\tilde l_2\circ_{\sigma}
\tilde l_2)(\alpha_1,\alpha_2,\alpha_3):=}
\\
&:=&
\sum_{\sigma\in \unsh(2,1)}e(\sigma)\epsilon(\sigma)\tilde
l_2(\tilde
l_2(\alpha_{\sigma(1)},\alpha_{\sigma(2)}),\alpha_{\sigma(3)})
\in d_H\om{N-1,0}
\end{eqnarray*}
\end{itemize}
for all $\alpha_1,\alpha_2,\alpha_3\in
\om{N,0}$,
where $\unsh(k,p)$ is the set of permutations $\sigma$ satisfying
\[
\underbrace {\sigma(1)<\dots<\sigma(k)}_{\rm first\ \sigma\ hand}
\hspace{.5cm}{\rm and}\hspace{.5cm} \underbrace
{\sigma(k+1)<\dots<\sigma(k+p),}_{\rm second\ \sigma\ hand}
\]
$e(\sigma)$ is the Koszul sign and 
$\epsilon(\sigma)$ is the standard sign of a permutation.
The meaning of the last condition is that $\tilde l_2$ satisfies the
Jacobi identity up to a boundary, that is
\[
 \tilde l_2(\tilde
l_2(\alpha_1,\alpha_2),\alpha_3) +
 \tilde l_2(\tilde
l_2(\alpha_2,\alpha_3),\alpha_1) +
\tilde l_2(\tilde
l_2(\alpha_3,\alpha_1),\alpha_2) \in d_H\om{N-1,0}.
\] 
The paper~\cite{barnich-fulp-lada-stasheff:preprint}
proves that when the original  map
$\tilde l_2$ satisfies conditions (i), (ii) and (iii),
there is a {\em strong homotopy Lie structure} 
(see~\cite{lada-markl:CommAlg95}
for relevant definitions)  on the graded
vector space $X_*=\Om*$, that is,
 a collection of linear, skew symmetric maps
$l_k:(X^{\otimes k})_*\longrightarrow X_{*+k-2}$, $k\geq 1$,
that satisfy, for any $n\geq 1$, the relation
\[
 \sum_\doubless{i+j=n+1}{\sigma\in \unsh(i,n-i)}
e(\sigma)\epsilon(\sigma)(-1)^{i(j-1)}
l_j(l_i(x_{\sigma(1)},\dots,
x_{\sigma(i)}),x_{\sigma(i+1)},
\dots ,x_{\sigma(n)}) =0.\nonumber
\]
In this formulation the element $l_1$  of degree $-1$ is the
boundary
operator of the complex. In the case under consideration,
$l_1=d_H.$

The authors  use the exactness of the complex $\Om*$
to define the lifting
$l_2$ as well as to  define the higher brackets, $l_k$. From
our point of view there is a problem with their construction
in that the recursive definition of $l_k$ is based on a choice
for each $k$-tuple of local functions, and there is
no control of the class of operator being defined.
For example, they argue:

{\em ``In degree zero, $\sum_{\sigma\in  \unsh(2,1)}
e(\sigma)\epsilon(\sigma)
\tilde l_2(\tilde
l_2(x_{\sigma(1)},x_{\sigma(2)}), x_{\sigma(3)})$
is equal to a boundary $b$ in $X_0$ by condition
(iii). There exists an element $z\in X_1$ with $l_1z=b$ and so we
define $\tilde l_3(x_1, x_2, x_3)=-z$.''}

This ``pointwise" construction is not sufficient to
guarantee if $\tilde l_2$ is a \ldo\
that the lift $\tilde l_3$ will be also.
Given a higher bracket $l_k$
which has been partially defined up to degree $j$
$$
l_k:\bigoplus_{i\leq j} (X^{\otimes k})_i\longrightarrow
\bigoplus_{i\leq j}X_{i+k-2},
$$
similar arguments are invoked to extend it to
$$
l_k:(X^{\otimes k})_{j+1}\longrightarrow X_{j+k-1}.
$$
There is no guarantee that the map so defined will
satisfy $l_k\in\dend_{k-2}(k)$.

The results of Sections 4 and 5 together with
Proposition~\ref{dtochar} given below allow
us to prove that in fact all $l_k\in\dend_{k-2}(k)$. To be more
precise,

\begin{itemize}
\item[--]
Theorem 4.8(ii) will establish
the existence of a suitable extension of $\tilde l_2$,
\item[--]
Theorem 5.2  will allow us to define $\tilde l_3\in O^{N-1}(3)$ as
well as the extension to $l_3\in \dend_1 (3)$.
\item[--]
Once we have defined $l_3$
the existence of the higher $l_k$ follows from Theorem 4.8 $(i)$,
i.e., the
acyclicity in dimensions greater than zero of $\dend_*(n)$.
\end{itemize}

\noindent
Recall that the {\em Euler operator\/} $E : \om{N,0} \to \loc(E)$
is defined by
\[
E(L(\x,\u)dx^1 \wedge \ldots \wedge dx^N) :=
\sum_I (-1)^I \left( \frac{d}{d\x}
\right)^I\left(\frac{\pa L}{\pa u_I}
\right).
\]
The next well-known lemma and the subsequent proposition play
an essential r\^ole here by allowing us to express in terms of
the character of a LDO conditions such as (i), (ii) and (iii)
which involve the image of $d_H$, that is, $d_H$ of
some unspecified forms, for a proof see \cite{Olver}.

\begin{lemma}
\label{Avro_Lancaster}
An element $\alpha\in\om{N,0}$ has the form $d_H\beta$ for
$\beta\in \om{N-1,0}$
if and only if $E(\alpha)=0.$
\end{lemma}

\begin{proposition}
\label{dtochar}
Given $A\in O^N(n)$,
$A(f_1,\ldots, f_n)\in d_H(\om{N-1,0})$ for all n-tuples
$(f_1,\ldots,f_n)\in Loc(E)^{\ot n}$
if and only if  $\chi(A)=0$.
\end{proposition}

\noindent
{\bf Proof. } From (\ref{again})
we know that there exists $B\in O^{N-1}(n)$
such that $A=d_H\circ B+\chi (A)$.
Since $E\circ d_H=0$ we have
\begin{equation}
\label{Avro_Shackelton}
E\circ \chi(A)=E\circ A.
\end{equation}
By this equation, $\chi(A)=0$ implies $E \circ A=0$ and
$A(f_1,\ldots, f_n)\in d_H(\om{N-1,0})$ by
Lemma~\ref{Avro_Lancaster}.

To prove the opposite implication, it is enough,
again by~(\ref{Avro_Shackelton}), to
show that $E\circ \chi(A)=0$
implies $\chi(A)=0$. Since $\chi(A)\in \ldo^0(n)$, this will follow
from the following lemma.

\begin{lemma}
\label{Supermarine_Spitfire}
For $A\in \ldo^0(n)$, $E\circ A = 0$ if and only if $A=0$.
\end{lemma}

\noindent
{\bf Proof.}
Let us discuss the linear case, $n=1$, first. Each exponent $\alpha$
can be written as $\alpha = (n,\beta)$, where $n:= \alpha(\rada
00)$ and $\beta$ is the rest of the array. 
Clearly $\deg_f(\alpha) = \deg_f(\beta)$. Thus each $A \in
\ldo^0(1)$ (which contain no horizontal derivatives)
decomposes to the sum $A = \sum_{m\geq 0} A_m$ with
\[
A_m := \sum_{n\geq 0,\ \deg_f(\beta)=m} 
p_{(n,\beta)}(\x,\u) \frac{\partial^n}{\partial u^n} 
 \left(\frac{\partial}{\partial \u} \right)^{\beta}.
\]
If $A \not = 0$, then the minimum $M := \min\{m;\
A_{m}\not=0\}$ is defined and the supremum 
\[
n_M := \sup\{n;\mbox{ there exists $\beta$ with $\deg_f(\beta)=M$ such
that }
p_{n,\beta}(\x,\u) \not=0\}
\] 
is, by~(\ref{dfsd}), a finite number.

It is immediate to see
that in the similar decomposition $E\circ A = \sum_{m\geq 0}(E\circ
A)_m$ one has $(E\circ A)_m = 0$ for $m < M$ and
\[
(E\circ A)_M  =
\sum_{\deg_f(\beta)=m} 
p_{(n_M,\beta)}(\x,\u) \frac{\partial^{n_M+1}}{\partial u^{n_M+1}} 
\left(\frac{\partial}{\partial \u} \right)^{\beta} +
\mbox{
terms of degrees $\leq n_M$ in ${\partial}/{\partial u}$.
}
\]
The assumption $E \circ A = 0$ implies that  $(E\circ A)_M = 0$, which
in turn implies that 
\[
\sum_{\deg_f(\beta)=M}p_{(n_M,\beta)}
\left(\frac{\partial}{\partial \u} \right)^{\beta}
=0,
\] 
which contradicts the definition of $n_M$.

Now in the multilinear case  the operators
in $A$ which are acting on $f_1$ look like the terms in the
case $n=1$,
involving only $(\fr\pa{\pa u^1})^\alpha$ and no total derivatives.
Essentially from the same argument
isolating the terms in $E\circ A$ which have
$\fr \pa{\pa u^1}$ we conclude that $A=0$.%
\qed

Proposition~\ref{dtochar} has the following corollary.

\begin{corollary}[\rm `pointwise = global']
\label{Opulka_se_Zebrulkou}
Given an element $A \in O^N(n)$, 
then $A(\Rada f1n) \in d_H\om{N-1,0}$ for each
$\Rada f1n \in \loc(E)$ if and only if $A$ is a boundary in $O^*(n)$.
\end{corollary}

Applying the corollary, we see that (ii) implies the existence of an
$\tilde m \in O^{N-1}(2)$ such that 
\[
\tilde l_2 +\tilde l_2\circ \tau=d_H\tilde m
\mbox { ($\tau$ interchanges the arguments).}
\]
We may assume that the coefficients of $\tilde m$ are symmetric as
bilinear operators, otherwise we replace $\tilde m$ by $\frac
12(\tilde m + \tilde m \circ \tau)$.
Replacing $\tilde l_2$ with $\tilde l_2 -\frac 12 d_H\tilde
m$ gives
a skew-symmetric element of $O^N(2)$ which determines the
same bracket on the space of functionals ${\cal F}$. Similar
arguments can be applied to insure appropriate skew-symmetries for
all other $l_k$'s. 
The next step is lifting $\tilde l_2$ to a chain
map on the
entire complex $X_*$.
Reinterpreting condition (i) with the
help of Proposition 6.4, we conclude that
\[
\chi(\tilde l_2\circ\fr d{dx^i_j})=0,
\]
for $i=1,2$ and $1\leq j \leq N.$
In the notation of Theorem \ref{Mixture}, $a_0(\tilde l_2)=0$
and  there exists a cycle  $l_2\in Z_0(\dend_*(2))$ such
that $\B_0(l_2)=\tilde l_2,$ that is,  a chain map of
LDO's lifting $\tilde l_2$.

Now ${\cal B}_0(\sum_{\sigma\in \unsh(2,1)}
l_2\circ_{\sigma} l_2) = \sum_{\sigma\in \unsh(2,1)}
\tilde l_2\circ_{\sigma}\tilde l_2$ and
condition (iii) implies that
\[
\chi(\sum_{\sigma\in \unsh(2,1)}
\tilde l_2\circ_{\sigma}\tilde l_2)=0.
\]
The existence of $l_3 \in \dend_1(3)$ now 
follows from Theorem~\ref{Mixture}(ii). 
We complete the strong homotopy Lie structure using the acyclicity
of $(\dend_*(n),\delta)$ in positive dimensions for all $n\geq 4.$

These results are summarized in the following theorem.
\begin{theorem}
\label{BFLS}
For $\alpha\in\om{N,0}$, let $\int \alpha$ be the functional
 $\int \alpha[\phi]=\int (j^\infty\phi)^* \alpha$.
Given  $\tilde l_2: \om{N,0}\otimes \om{N,0} \rightarrow \om{N,0},$
define a bilinear map
$\{-,-\} : {\cal F} \otimes {\cal F} \to {\cal F}$
on the space of functionals by
\begin{equation}
\label{bra}
\{\int \alpha_1,\int \alpha_2\}:=
\int \tilde l_2(\alpha_1,\alpha_2).
\end{equation}
If  $\tilde l_2$ is a LDO
satisfying conditions (i), (ii) and (iii)
listed at the beginning of this section,
then~(\ref{bra}) defines a Lie bracket which
lifts to  a strong homotopy Lie algebra structure on
$\om{*,0}$ such that all the higher brackets are also LDOs.
\end{theorem}
In other words, all the constructions in
\cite{barnich-fulp-lada-stasheff:preprint} can in fact be
done in the class LDO of local differential operators.
\nocite{BV1,bv:anti,bv:antired,dickey:ex}

%\bibliography{b}

\vskip3mm
\catcode`\@=11
\noindent
M.~M.: Mathematical Institute of the Academy, \v Zitn\'a 25, 115 67
Praha 1, Czech Republic,\hfill\break\noindent
\hphantom{M.~M.:} email: {\tt
markl@math.cas.cz}\hfill\break\noindent

\noindent
S.~S.: Department of Mathematics, University of Bar Ilan, Israel,
\hfill\break\noindent
\hphantom{S.~S.:} email: {\tt shnider@bimacs.cs.biu.ac.il}

\vfill

\hfill {\tt \jobname.tex}
\end{document}